\documentclass[12pt]{amsart}
\usepackage{latexsym,enumerate}
\usepackage{amsmath,amsthm,amsopn,amstext,amscd,amsfonts,amssymb,mathrsfs}
\usepackage{color}
\usepackage{color}
\usepackage{graphicx}

\numberwithin{equation}{section}

\newtheorem{theorem}{Theorem}
\newtheorem{lemma}{Lemma}

\newtheorem{remark}{Remark}

\numberwithin{theorem}{section}
\numberwithin{corollary}{section}
\numberwithin{lemma}{section}
\numberwithin{definition}{section}
\numberwithin{proposition}{section}
\numberwithin{remark}{section}

\newcommand*\rup{\mathbb{R}_+ ^2 }
\newcommand*\re{\mathbb{R}}

\newcommand*\N{\mathbb{N}}

\date{}

\title[Some isoperimetric inequalities with respect to monomial weights]{Some isoperimetric inequalities\\with respect to monomial weights}

\begin{document}

\author[A. Alvino]{A. Alvino$^1$}
\author[F. Brock]{F. Brock$^2$}
\author[F. Chiacchio]{F. Chiacchio$^1$}
\author[A. Mercaldo]{A. Mercaldo$^1$}
\author[M.R. Posteraro]{M.R. Posteraro$^1$}


\setcounter{footnote}{1}
\footnotetext{Universit\`a di Napoli Federico II, Dipartimento di Matematica e Applicazioni ``R. Caccioppoli'',
Complesso Monte S. Angelo, via Cintia, 80126 Napoli, Italy;\\
e-mail: {\tt angelo.alvino@unina.it, fchiacch@unina.it,  mercaldo@unina.it, posterar@unina.it}}

\setcounter{footnote}{2}
\footnotetext{Department of Mathematics, Computational Foundry, College of Science, Swansea University, Bay Campus, Fabian Way, Swansea SA1 8EN, Wales, UK, 
e-mail: {\tt friedemann.brock@swansea.ac.uk}
}

\begin{abstract} 
We solve a class of isoperimetric problems on
 $\mathbb{R}^2_+  :=\left\{ (x,y)\in \mathbb{R} ^2 :   y>0  \right\}$ with respect to monomial weights. Let $\alpha $ and $\beta $ be real numbers such that $0\le \alpha <\beta+1$, $\beta\le 2 \alpha$. 
We show that,  among all smooth sets $\Omega$ in $\mathbb{R} ^2_+$ with fixed weighted measure $\iint_{\Omega } y^{\beta} dxdy$, the weighted perimeter $\int_{\partial \Omega } y^\alpha \, ds$ achieves its minimum for a smooth set which is symmetric w.r.t. to the $y$--axis, and is explicitly given. Our results also imply an estimate of a weighted Cheeger constant and a lower bound for the first eigenvalue of a class of nonlinear problems.

\medskip

\rm 
\\[0.1cm]
\rm 
\end{abstract}
\maketitle

\vspace{2cm}

\section{Introduction}
The last two decades have seen a growing interest in isoperimetric inequalities with respect to weights.
\\
In most cases, volume and perimeter in those inequalities carried the same weight, because such a setting corresponds to {\sl manifolds with density}. However, most research dealt with inequalities where both the volume functional and perimeter functional carry the same weight, see  for instance \cite{BCMR}, 
\cite{BBMP2}, \cite{BBCLT}, \cite{BMP}, \cite{DFP},   \cite{XR},   \cite{XRS}, \cite{CMV},   \cite{CJQW}, \cite{Cham}, 
\cite{DDNT}, 
\cite{KZ}, \cite{MadSalsa}, \cite{Mo}, \cite{Mo2}, \cite{MP}  and  the  references therein.\\
More recently, also problems with  {\sl different} weight functions for perimeter and volume were studied, see for example \cite{ABCMP}, \cite{ABCMP_atti}, \cite{ABCMP_die},  \cite{BBMP},    \cite{C},   \cite{diGiosia_etal},  \cite{DHHT}, \cite{Howe}, \cite{M}, \cite{MS}, \cite{PS} and  the  references therein. However, there is only  a sparse literature on situations where the isoperimetric sets are not radial, see \cite{DHHT}, \cite{Castro}, \cite{abreu}.

In this paper we study the following isoperimetric problem:
\bigskip 

{\sl Minimize $\displaystyle\int_{\partial \Omega } y^\alpha ds $ 
among all smooth sets 
$\Omega \subset \mathbb{R} ^{2}_{+}$ satisfying $\displaystyle\iint_\Omega
 y^\beta \, dxdy =1$}
\medskip

or equivalently
$$
{\bf (P)} \qquad \qquad \inf 
\left\{ \displaystyle{ 
\frac{ \displaystyle\int_{\partial \Omega } y^\alpha ds}
{ \displaystyle{ \left[ \iint_\Omega  y^\beta \, dxdy \right] ^{(\alpha +1)/(\beta +2 )}}}}: \,
\, 0<\iint_\Omega  y^\beta \, dxdy <+\infty \right\} 
=: \mu (\alpha , \beta ).
$$  
Our main result, proved in Section 2, is the following:

\begin{theorem}\label{maintheorem}
Assume that 
\begin{equation}
\label{cond11}
  0  \leq \alpha <\beta +1 
\end{equation}
and
\begin{equation}
\label{cond21}
 \beta \leq 2\alpha.
\end{equation}
Then problem {\bf (P)} has a minimizer which is given by 
\begin{eqnarray}
\label{minU}
 \Omega^\star  & := & \{ (x,y) :\, |x|< f (y) , \, 0 < y <1  \},
\quad \mbox{
where }
\\
\label{xi}
f (y) & := & \int_y ^1 
\frac{t^{\beta -\alpha +1}}{\sqrt{1-t^{2(\beta -\alpha +1)} }} \, dt , \quad (0<y<1).
\end{eqnarray}
Moreover, we have
\begin{equation}
\label{muab}
\mu (\alpha ,\beta )=
\gamma ^{ 
\frac{\alpha +1}{\beta +2} -1 
} 
\cdot
\left[
\frac{
(\beta +1)(\beta +2)
}{
\alpha +1 
} 
\right] ^{\frac{\alpha +1}{\beta +2}} 
\cdot 
\left[ 
B 
\left( 
\frac{\alpha +1}{2\gamma } , \frac{1}{2} 
\right) 
\right] ^{\frac{\gamma}{\beta +2 }} ,
\end{equation}
where $\gamma := \beta +1 -\alpha $ and $B$ denotes the Beta function.
In particular,
\begin{equation}
\label{mualpha2alpha}
\mu (\alpha , 2\alpha ) = \sqrt{\frac{2\pi (2\alpha +1)}{\alpha +1}}.
\end{equation}
\end{theorem}

\begin{remark}\label{Remark 1.1}\rm  

\noindent {\bf (a)} 
First observe  that $\Omega^\star$ is the half-circle when $\alpha=\beta$. Therefore  Theorem \ref{maintheorem} includes the result obtained by Maderna and Salsa in \cite{MadSalsa} (see also \cite{XR}, \cite{BCM3}).
\\
{\bf (b)} Let 
$B_{1}^{+} := \{ (x,y)\in \mathbb{R} ^2 :\, x^2 + y^2 <1 , \,  y>0\} $. It is elementary to verify that,
\\
\noindent{\bf 1.} if $\beta -\alpha <0$ then $\Omega \supseteq B_{1}^{+}$, 
\\
\noindent{\bf 2.} if $\beta -\alpha >0$ then $\Omega \subseteq B_{1}^{+}$, 
\\
\noindent{\bf 3.} if  $\beta -\alpha =0$ then  $\Omega \equiv B_{1}^{+}$; see Figure 1.
\end{remark}
\medskip

\begin{figure}[h]
 \includegraphics[scale=0.445]{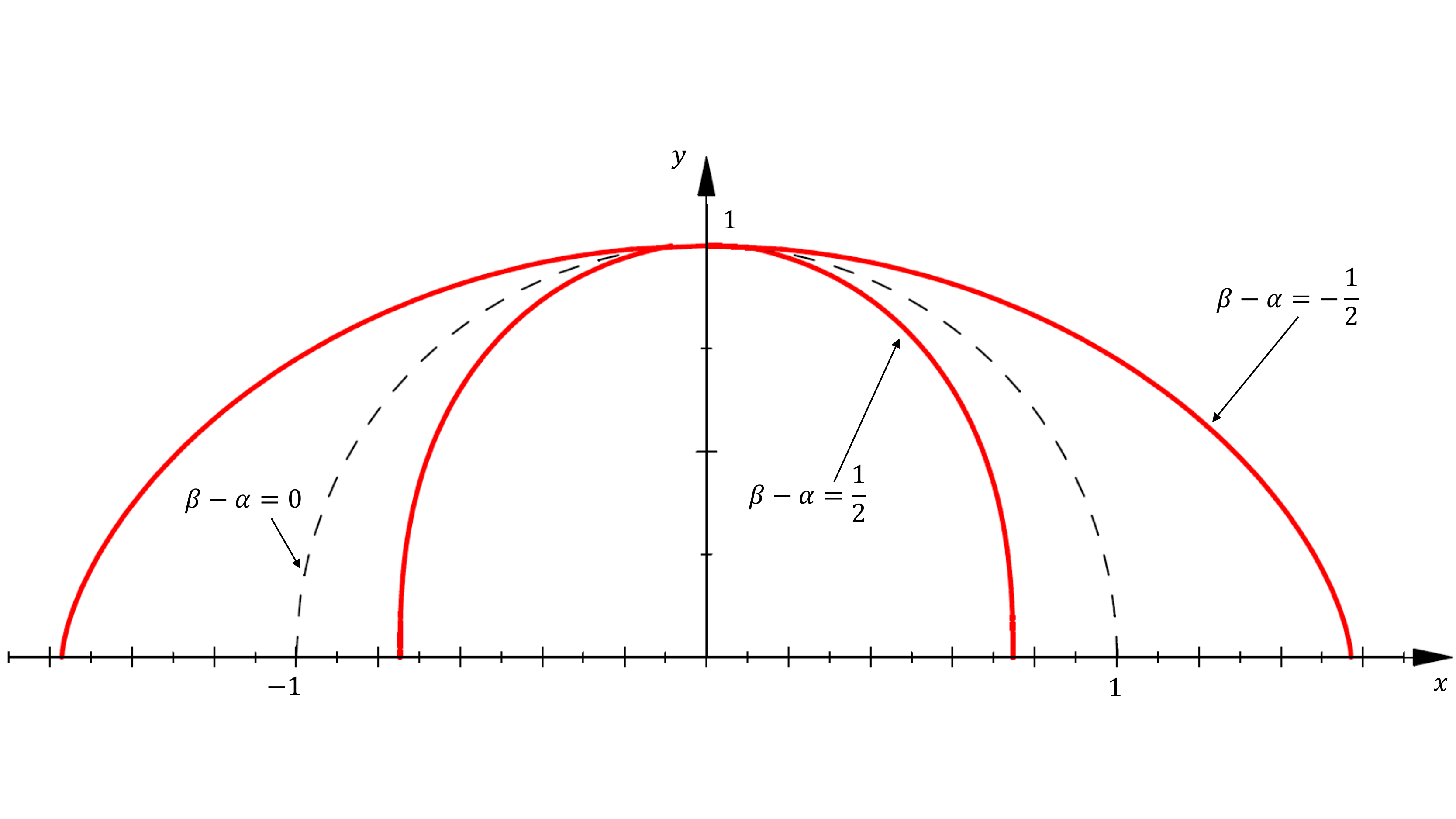}
\caption{ Isoperimetric sets for different values of $\alpha $ and $\beta $}
\end{figure}

Theorem 1.1  also allows to obtain a {\sl  Faber-Krahn} - type inequality for the so-called  {\sl weighted Cheeger constant}, and in turn a lower bound for the first eigenvalue  for a degenerate elliptic operators. For similar results see also \cite{BDNT}, \cite{BF1}, \cite{BF2}, \cite{CMN},  \cite{IL}, \cite{P}, \cite{Saracco}.

\section{Isoperimetric inequality in the upper half plane}
Let $\rup := \{ (x,y)\in \mathbb{R} ^2 :\, y>0\} $. Throughout this paper, we assume that $\alpha ,\beta \in \mathbb{R} $ and
\begin{equation}
\label{alphabeta}
\beta +1>0  \ \mbox{ and } \ \alpha \geq 0.
\end{equation} 
If 
$\Omega \subset \rup $ is measurable,  we set
\begin{eqnarray}\label{def}
\Omega (y) & := & \{ x\in \mathbb{R}:\, (x,y)\in \Omega \} , \quad (y\in \mathbb{R}_+  ),
\\
\Omega ^{\prime } & := & \{ y\in \mathbb{R} _+ :\, \Omega (y)\not= \emptyset \} .
\end{eqnarray}
Further, we define the {\sl weighted area } of $\Omega $ by 
$$
A_{\beta } (\Omega ) :=
\iint_\Omega y^{\beta } \, dxdy ,
$$
and the {\sl weighted relative perimeter} of $\Omega $ as
\begin{equation}
P_{\alpha }(\Omega  ) :=
\sup \left\{ \iint_{\rup} \mbox{div}\,\left( y^{\alpha }
\mathbf{v}\right) \,dx dy:\,\mathbf{v}\in C_{0}^{1}(\rup ,\mathbb{R}^2),\,|
\mathbf{v}|\leq 1\mbox{ in }\, \Omega \right\} .
\end{equation}
It is well-known that, if $\Omega $ is an open, rectifiable set, then the  following equality holds
\begin{equation}
P_{\alpha }(\Omega )= 
\int_{\partial \Omega \cap \rup } y^{\alpha } \, d\mathscr{H} _{1}
\ .
\end{equation}
($\mathscr{H} _1 $ denotes $1$-dimensional Hausdorff-measure.)
\\[0.1cm]
\begin{remark}\rm The following properties of the perimeter are well-known:
\\
Let $\Omega $ be measurable with  $0<A_{\beta }(\Omega )<+\infty $ and $P_{\alpha } (\Omega) <+\infty$. 
\\
Then there exists a sequence of open, rectifiable sets $\{ \Omega _n \} $ with $\displaystyle\lim_{n\to \infty } A_{\beta }(\Omega \Delta \Omega _n )=0$ and
\begin{equation}
\label{limperimeter}
P_{\alpha } (\Omega )= \lim _{n\to \infty } P_{\alpha } (\Omega _n ).
\end{equation}
Further, we have
\begin{equation}
\label{lscperimeter}
P_{\alpha } (\Omega )\leq  \liminf_{n\to \infty } P_{\alpha } (\Omega _n )
\end{equation} 
for any sequence of open, rectifiable sets $\{ \Omega _n \}$ satisfying  $\displaystyle\lim_{n\to \infty } A_{\beta }(\Omega \Delta \Omega _n )=0$.
\end{remark} 
\bigskip
    
We define the ratio
$$
{\mathcal R}_{\alpha , \beta } (\Omega )  := \displaystyle{ 
\frac{ P_{\alpha }( \Omega  ) }{ \displaystyle{ \left[ A_{\beta }(\Omega ) \right] ^{(\alpha +1)/(\beta +2 )}}}}, \quad (0<A_{\beta } (\Omega ) <+\infty ).
$$
\begin{remark}\label{remark2.2}\rm  We have ${\mathcal R}_{\alpha , \beta } (t\Omega ) = {\mathcal R} _{\alpha , \beta } (\Omega )$ for every $t>0$.
\end{remark}
\bigskip

We study the following isoperimetric problem:
$$
{\bf (P)} \qquad \qquad \inf \left\{ {\mathcal R}_{\alpha , \beta }  (\Omega ): \, 0<A_\beta(\Omega ) <+\infty \right\} =: \mu (\alpha , \beta ).
$$  
Our first aim is to reduce the class of admissible sets in the isoperimetric problem {\bf (P)}.
\\
Throughout our proofs let $C$ denote a generic constant which may vary from line but does not depend on the other parameters.
\\ 
The first two Lemmata give necessary conditions for a minimizer to exist.
\medskip

\begin{lemma} \label{lemma2.3} If $\alpha >\beta +1$, then 
\begin{equation}
\label{inf01}
\mu (\alpha , \beta )= 0, 
\end{equation}
and {\bf (P)} has no minimizer.
\end{lemma}
\noindent {\sl Proof: } 
Let $\Omega (t) := (0, t) \times (0,1)$, ($t>0$). Then
\begin{eqnarray*}
\int_{\partial \Omega(t)\cap \rup } y^{\alpha }\, ds & = & t+2 \int_0 ^t y^{\alpha} \, dy = t+ \frac{2}{\alpha +1} ,
\\
 \iint\limits_{\Omega (t) } y^{\beta }\, dxdy & = &
\frac{t}{\beta +1 }.
\end{eqnarray*}
Hence
$$
{\mathcal R} _{\alpha , \beta } (\Omega (t)) = 
\frac{t+ (2/(\alpha +1))}{[t/(\beta +1)]^{(\alpha +1)/(\beta +2)} }\longrightarrow 0, \mbox{ as $t\to + \infty $,}
$$
and the assertion follows.
$\hfill \Box $

\begin{lemma}\label{lemma2.4}  Assume $\alpha <\beta +1   $. Further, let $\Omega \subset \mathbb{R}^2 _+ $ be a nonempty, open and rectifiable set, which is not simply connected. Then there exists a nonempty, open and  rectifiable set $U \subset \mathbb{R}^2 _+ $ which is simply connected, such that
\begin{equation}
\label{UOmega} 
{\mathcal R} _{\alpha , \beta } (U) < {\mathcal R} _{\alpha , \beta }(\Omega ).
\end{equation}
\end{lemma}
\noindent {\sl Proof: } {\bf (i)} First assume that $\Omega $ is connected.  Let $G $ be the unbounded component of $\mathbb{R} ^2 \setminus \overline{\Omega  }$ and set $U := \mathbb{R} ^2 \setminus \overline{G }$. Then $U $ is simply connected with
$\Omega \subset U \subset \mathbb{R} _+ ^2 $ and $\partial U \subset \partial \Omega  $, so that (\ref{UOmega}) follows.
\\
{\bf (ii)}
Next, let
$\Omega = \cup_{k=1} ^m \Omega _k $, with mutually disjoint, nonempty, open, connected and rectifiable sets $\Omega _k $, ($k=1,\ldots ,m$, $m\geq 2$). We set ${\mathcal R} _{\alpha , \beta } (\Omega )=: \lambda $. Let us assume that ${\mathcal R} _{\alpha , \beta } (\Omega _k )\geq \lambda $ for every $k\in \{ 1, \ldots ,m\} $. Then we have, since $(\alpha +1 )/(\beta +2 )<1$,
\begin{eqnarray*}
P_{\alpha }(\Omega  ) & = & \sum _{ k=1 }^m P_{\alpha } (\Omega _k  ) \geq  \lambda \left[ \sum_{k=1} ^m A_{\beta } (\Omega _k )\right] ^{(\alpha +1 )/(\beta +2)} 
\\
 & > & \lambda \left[ \sum_{k=1} ^m A_{\beta } (\Omega _k ) \right] ^{(\alpha +1 )/(\beta +2)}  = \lambda \left[ A_{\beta } (\Omega )\right] ^{(\alpha +1 )/(\beta +2)} 
\\
 & = & P_{\alpha }(\Omega  ),
\end{eqnarray*}
a contradiction. Hence there exists a number $k_0 \in \{1, \ldots , m\} $ with ${\mathcal R}_{\alpha , \beta } (\Omega _{k_0 }) < \lambda $. Then, repeating the argument of part {\bf (i)}, with $\Omega _{k_0 }$ in place of $\Omega $, we again arrive at (\ref{UOmega}).
$\hfill \Box $
\begin{lemma}\label{lemma2.5} There holds  
\begin{equation}
\label{betaplus1}
\mu (\beta +1 , \beta ) = \beta +1 ,
\end{equation}
but {\bf (P)} has no open rectifiable minimizer.
\end{lemma}

\noindent{\sl Proof: } With $\Omega (t)$ as in the proof of Lemma \ref{lemma2.3}, we calculate
${\mathcal R} _{\beta +1 , \beta }(\Omega (t)) >\beta +1$ and 
\begin{equation}
\label{limR1} 
\lim_{t\to \infty } {\mathcal R} _{\beta +1 , \beta } (\Omega (t) ) = \beta +1.
\end{equation}
Let $\Omega \subset \mathbb{R}_+ ^2 $ be open, rectifiable and simply connected. Then $\partial \Omega $ is a closed Jordan curve ${\mathcal C}$ 
with counter-clockwise representation
$$
{\mathcal C} = \{ (\xi (t), \eta (t)): \, 0\leq t\leq a \}, \ \ (a\in 0, +\infty ),
$$
where $\xi , \eta \in C[0,a] \cap C^{0,1} (0,a)$, $\xi (0)= \xi (a)$, $\eta (0)= \eta (a)$, and $(\xi ') ^2 + (\eta ')^2 >0 $ on $[0,a]$. Using Green's Theorem we evaluate
\begin{eqnarray}
\nonumber
A_{\beta }(\Omega )  & = & - \int_{{\mathcal C} } \frac{y^{\beta +1 }}{\beta +1 } \, dx = -\int_0 ^a \frac{[\eta (t) ] ^{\beta +1}}{\beta +1}  \xi ' (t)\, dt 
\\
\label{ineqconnected}
 & \leq & \frac{1}{\beta +1 } \int_0 ^a  [\eta (t) ] ^{\beta +1} \sqrt{ (\xi ' )^2 + (\eta ')^2 } \, dt = \frac{1}{\beta +1} P_{\beta +1 }(\Omega  ).  
\end{eqnarray}   
Equality in 
(\ref{ineqconnected}) can hold only if 
$
\eta '   \equiv 0 $ and $\xi '  \leq 0 $ on $\partial \Omega \cap \mathbb{R} ^2_+ $, that is, if  $\partial \Omega \cap \mathbb{R} ^2_+ $ is a single straight segment which is parallel to the $x$-axis. But this is impossible. Hence we find that
\begin{equation}
\label{strictineqisop}
P_{\beta +1 }(\Omega  ) > (\beta +1) A_{\beta }(\Omega ).
\end{equation}
To show the assertion in the general case, we proceed similarly as in the proof of Lemma \ref{lemma2.4}:
\\ 
Assume first that $\Omega $ is connected and define the sets $G$ and $U$ as in the last proof.  Using 
(\ref{strictineqisop}), with $U$ in place of $\Omega $, we obtain
\begin{eqnarray}
\nonumber
P_{\beta +1 }(\Omega  ) & \geq & 
P_{\beta +1 }(U  ) 
\\
\label{strict2}
 & > & (\beta +1 )A_{\beta } (U )  \geq  (\beta +1) A_{\beta }(\Omega  ).
\end{eqnarray} 
Finally, let $\Omega $ be open and rectifiable. Then $\Omega = \cup_{k=1 } ^m \Omega _k $, with mutually disjoint, connected sets $\Omega _k $, ($k=1, \ldots ,m$). Then (\ref{strict2}) yields
\begin{eqnarray}
\nonumber
P_{\beta +1 } (\Omega  ) & = & \sum_{k=1 } ^m P_{\beta +1 }(\Omega _k  ) 
\\
\label{strict3}
 & > &
(\beta +1) \sum_{k=1}^m A_{\beta } (\Omega _k ) = (\beta +1 )A_{\beta }(\Omega ).    
\end{eqnarray} 
Now the assertion follows from (\ref{strict3}) and (\ref{limR1}).
$\hfill \Box $

\begin{lemma}\label{lemma 2.6}  
Let $ 2\alpha <\beta $. Then (\ref{inf01}) holds and {\bf (P)} has no minimizer.
\end{lemma}
{\sl Proof: } 
Let $z(t):= (0, t)$, ($t\geq 2$). 
Then  we have for all $t\geq 2$, 
\begin{eqnarray*}
\int_{\partial B_1 (z(t)) } y^{\alpha } ds & \leq & 
C t^{\alpha }  \ \mbox{ and }
\\
\iint\limits_{B_1 (z(t)) }  y^{\beta } \, dxdy  & \geq & C t^{\beta } .
\end{eqnarray*}
This implies
\begin{equation}
\label{lim02}
{\mathcal R}_{\alpha , \beta } (B_1 (z(t))) \leq  C t^{\alpha - \beta (1+\alpha )/(2+\beta )} \longrightarrow 0, \ \mbox{ as }\ t\to +\infty ,
\end{equation}
and the assertion follows.
$\hfill \Box $

\bigskip

Next we recall the definition of the Steiner symmetrization w.r.t. the $x$-variable.
If $\Omega $ is measurable, we set
$$
S(\Omega ) := \{ (x,y) :\, x\in S(\Omega (y)), \, y\in \Omega ^{\prime} \}
$$
where $\Omega (y)$, $\Omega ^{\prime}$ are defined in \eqref{def} and 
$$
S(\Omega (y) ) :=
\left\{ 
\begin{array}{ll}
\left( -\frac{1}{2}\mathscr{L}^1 (\Omega (y)) , + \frac{1}{2}\mathscr{L}^1 (\Omega (y)) \right) & \mbox{ if }\ 0<\mathscr{L}^1 (\Omega (y)) <+\infty 
\\
\emptyset & \mbox{ if }\ \mathscr{L}^1 (\Omega (y)) =0
\\
\mathbb{R} & \mbox{ if } \  \mathscr{L}^1 (\Omega (y)) +\infty 
\end{array}
\right.
.
$$
Note that $S(\Omega) (y)$ is a symmetric interval with 
$ \mathscr{L}^1 (S(\Omega )(y)) =\mathscr{L}^1 (\Omega (y))$.
\\
Since the weight functions in the functionals $P_{\alpha }$ and $A_{\beta }$ do not depend on $x$, we  have the following well-known properties, see \cite{HarHoweMor}, Proposition 3. 

\begin{lemma}
\label{lemma 2.7} Let $\Omega  \subset \rup $ be measurable.
Then
\begin{eqnarray}
\label{perimeterrearrangeS}
 & & P_{\alpha } (\Omega ) \geq P_{\alpha }(S(\Omega ) ) \quad \mbox{and}
\\
\label{arearearrangeS}
 & & A_{\beta }(\Omega ) = A_{\beta }(S(\Omega )).
\end{eqnarray}  
For nonempty open sets $\Omega $ with $\Omega = S(\Omega )$ we set
$\Omega _t := \{ (x,y)\in \Omega :\, y>t\} $ and 
\begin{eqnarray}
\label{y+}
y^+ & := & \inf \{ t\geq 0 :\, A_{\beta }( \Omega _t ) =0 \},
\\
\label{y-}
y^- & := & \sup \{ t\geq 0 :\, A_{\beta }( \Omega \setminus \Omega _t ) =0 \}
.
\end{eqnarray} 
\end{lemma}

\noindent
\\
{\bf Remark 2.3.}
Assume that $\Omega \subset \rup $ is a bounded, open and rectifiable set with
$0<A_{\beta }(\Omega ) <+\infty $ and $\Omega = S(\Omega )$.
Then it has the following representation,
\begin{eqnarray}
\label{repres1}
\Omega & = & \{ (x,y):\, |x|<f(y), \, y^- <y<y^+ \}, \quad \mbox{where  }
\\
\nonumber
 & & f: (y^-,y^+ ) \to (0, +\infty ) \ \mbox{ is lower semi-continuous.} 
\end{eqnarray}
\begin{lemma}
\label{lemma2.9}
Let $\Omega $ be a nonempty, bounded, open and rectifiable set with $\Omega = S(\Omega )$. Then we have
\begin{eqnarray} 
\label{est1}
P_{\alpha } (\Omega  ) & \geq & \frac{2}{\alpha +1}  ((y^+ )^{\alpha +1} -(y^-) ^{\alpha +1 }) \quad \mbox{ and}
 \\
\label{est2}
 P_{\alpha }(\Omega  ) & \geq & 2 y^{\alpha } f(y) \quad \forall y\in (y^- , y^+ ),
\end{eqnarray}
where $f$ is given by (\ref{repres1}).
\end{lemma}

\noindent {\sl Proof: } Assume first that $\Omega $ is represented by (\ref{repres1}) where
\begin{equation}
\label{fsmooth}
f\in C^1 [y^- , y^+ ] \ \mbox{ and } f(y^- )= f(y^+ )=0.
\end{equation}
Then we have for every $y\in (y^- , y^+ )$,
\begin{eqnarray*}
P_{\alpha } (\Omega ) & = & 2\int_{y^- }^{y^+ } t^{\alpha } \sqrt{1+ (f^{\prime }  (t))^2}\, dt  \geq  2\int_{y }^{y^+ } t^{\alpha } \sqrt{1+ (f^{\prime }  (t))^2}\, dt  
\\
 & \geq & 2 y^{\alpha } \int_{y} ^{y^+ } |f^{\prime } (t)|\, dt \geq 2 y^{\alpha } f(y).
\end{eqnarray*}
Furthermore, there holds
\begin{eqnarray*}
P_{\alpha } (\Omega ) & = & 2\int_{y^- }^{y^+ } t^{\alpha } \sqrt{1+ (f^{\prime }  (t))^2}\, dt  \geq  2\int_{y^-  }^{y^+ } t^{\alpha } \, dt  
\\
 & = & \frac{2}{\alpha +1}  ((y^+ )^{\alpha +1} -(y^-) ^{\alpha +1 }
) .
\end{eqnarray*}
In the general case the assertions follow from these calculations by approximation with sets $\Omega $ of the type given by  (\ref{repres1}), (\ref{fsmooth}).
$\hfill \Box $

\begin{lemma}
\label{lemma2.10}Assume that 
\begin{eqnarray}
\label{cond1}
 & & \alpha <\beta +1, \quad \mbox{and } 
\\
\label{cond2}
 & & \beta <2\alpha ,
\end{eqnarray} 
and let $\Omega $ be a bounded,  
open and rectifiable set with 
$A_{\beta }(\Omega )=1 $, $\Omega = S(\Omega )$ and $P_{\alpha } (\Omega )<\mu (\alpha , \beta )+1  $. 
Then there exist  positive numbers $C_1 $ and $C_2 $ which depend only on $\alpha $ and $\beta $ such that  
\begin{equation}
\label{boundy+}
C_1 \geq y^+  \quad \mbox{ and }\quad y^+ - y^- \geq C_2 .
\end{equation}
\end{lemma}

\noindent{\sl Proof: } By  (\ref{est1}) and (\ref{est2}) we have 
\begin{eqnarray}
\label{est3}
\frac{1}{2} (\mu (\alpha , \beta )+1) (\alpha +1) & \geq & (y^+)^{\alpha +1} - (y^-)^{\alpha +1} \ \mbox{ and }
\\
\label{est4}
\frac{\mu (\alpha , \beta )+1}{2} & \geq & y^{\alpha } f(y) \quad \forall y\in ( y^- , y^+ ).
\end{eqnarray}
It follows that
\begin{eqnarray}
\nonumber
1 & = & A_{\beta } (\Omega ) = 2 \int_{y^- } ^{y^+} y^{\beta } f(y)\, dy 
\\
\label{est5}
 & \leq & (\mu (\alpha , \beta )+1 ) \int_{y^- } ^{y^+ } y^{\beta -\alpha } \, dy = \frac{\mu (\alpha , \beta )+1 }{\beta +1-\alpha } ((y^+)^{\beta +1 -\alpha } -(y^- )^{\beta +1 -\alpha } ).  
\end{eqnarray}
Setting $z := \frac{y^- }{y^+}  (\in [0,1 ))$, we obtain from (\ref{est3}) and (\ref{est5}),
\begin{eqnarray}
\label{est6}
(y^+ )^{\alpha +1 } & \leq &
\frac{(\mu (\alpha , \beta ) +1) (\alpha +1)
}{ 
2(1-z^{\alpha +1 })
} \ \mbox{ and }
\\
\nonumber
(y^+ )^{\beta +1-\alpha } & \geq & \frac{
\beta +1-\alpha 
}{ 
(\mu (\alpha , \beta )+1) (1-z ^{\beta +1-\alpha })} ,
\end{eqnarray}
which implies that 
$$
f(z ):= \frac{(1-z ^{\alpha +1})^{(\beta +1-\alpha )/(\alpha +1)}}{ 1-z ^{\beta +1-\alpha }} \leq C, 
$$
with a constant $C$ which depends only on $\alpha $ and $\beta $.
By (\ref{cond2}) we have that 
$$
\lim _{z\to 1-} f(z) =+\infty .
$$
 Hence it follows that 
\begin{equation}
\label{est7}
z = \frac{y^- }{y^+ } \leq 1-\delta \ \mbox{ for some $\delta \in (0,1)$.}
\end{equation}
Using (\ref{est6}) and (\ref{est7}) this leads to (\ref{boundy+}).   
$\hfill \Box $ 

\begin{lemma}\label{lemma2.11} Assume (\ref{cond1}) and (\ref{cond2}).   
Then problem {\bf (P)} has a minimizer $\Omega^\star $ which is symmetric w.r.t. the $y$-axis. 
\end{lemma}

\noindent {\sl Proof: } 
We proceed in $4$ steps.
\\
{\sl Step 1: } {\sl A minimizing sequence :}
\\
Let $\{ \Omega _n \} $ a minimizing sequence, that is, $\lim _{n\to \infty } \mathcal{R} _{\alpha , \beta } (\Omega _n ) = \mu (\alpha , \beta )$.
 In view of the Remarks 2.1 and 2.2 and the Lemmata \ref{lemma2.4} and \ref{lemma2.9} we may assume that $\Omega _n $ is open, simply connected and rectifiable with $\Omega_n = S(\Omega _n )$, $A_{\beta } (\Omega _n )=1 $ and $P_{\alpha } (\Omega _n )\leq \mu (\alpha , \beta ) + \frac{1}{n} $, ($n\in \mathbb{N} $).
\\[0.1cm]
{\sl Step 2: } {\sl Parametrization of $\partial \Omega _n $ :} 
\bigskip

\noindent Denote  ${\mathcal C}_n :=\overline{\partial \Omega _n \cap \{ (x,y ):\, x>0 , \, y>0 \} }$. It is clear that ${\mathcal C}_n$ a simple smooth curve with
\begin{equation}
\label{arcl}
{\mathcal C}_n = \{ (x_n (s), y_n (s) ):\, s\in [0, L_n ]\},
\end{equation}
where $s$ denotes the usual arclength parameter, $L_n \in (0, +\infty)$, $x, y\in C[0, L_n ] \cap C^{0,1} (0, L_n )$,  $x_n (s) >0$ and $y_n (s) >0$ for every $s\in (0,L_n )$,   ($n\in \mathbb{N}$). We orientate ${\mathcal C}_n $ in such a way that the mapping $s\mapsto y_n (s)$ is nonincreasing and $x_n (0)=0 $. Setting $y_n (0) =: y_n ^+ $ and $y_n (L_n )=: y_n ^- $, we have by  Lemma \ref{lemma2.9},  
\begin{equation}
\label{ynbound} 
C_1 \geq y_n ^+ \quad \mbox{and }\ y_n ^+ - y_n ^- \geq C_2 ,
\end{equation}  
where $C_1 $ and $C_2 $ do not depend on $n $. Note  that $x_n (L_n )=0$ in case that $y_n ^- >0$.  Further, Lemma \ref{lemma2.9}, (\ref{est2}) shows that
\begin{equation}
\label{estxn}
\mu (\alpha , \beta )+\frac{1}{n} \geq 2 y_n(s)^{\alpha } x_n (s), \quad \forall s\in [0, L_n ).
\end{equation} 
For our purposes it will be convenient to work with another parametrization 
of ${\mathcal C}_n $:
We set 
\begin{eqnarray*}
\varphi _n (s) & := & 
\frac{2}{P_{\alpha } (\Omega  _n )}
\int _0 ^s y_n ^{\alpha } (t)\, dt  , \quad \mbox{and }
\\
X_n (\varphi_n (s)) & := & x_n (s), \ Y_n (\varphi _n (s)):= y_n (s), \quad (s\in [0, L_n] ).
\end{eqnarray*}
Then $X_n , Y_n \in C[0,1] \cap C^{0,1} (0,1)$, 
 and we evaluate
\begin{eqnarray}
\label{An} 
\quad
1  =  A_{\beta }( \Omega _n ) & = & \frac{2}{\beta +1 } \int_0 ^1 Y_n ^{\beta +1 } (\sigma ) X_n ^{\prime} (\sigma )\, d\sigma
=  - 2 \int_0 ^1 Y_n ^{\beta } (\sigma ) Y_n ^{\prime } (\sigma ) X_n (\sigma )\, d\sigma ,
\\
\label{Pn} 
P_{\alpha } (\Omega _n ) & = & 2\int_0 ^1 Y_n ^{\alpha } (\sigma ) \sqrt{ (X_n ^{\prime } (\sigma ))^2 + (Y_n ^{\prime} (\sigma ))^2 } \, d\sigma,
\\
\label{Xn}
\frac{d}{d\sigma } X_n (\sigma ) & = & \frac{P_{\alpha } (\Omega _n )}{2} (y_n (s))^{-\alpha } x_n ' (s),
\\
\label{Yn}   
\frac{d}{d\sigma } Y_n (\sigma ) & = & \frac{P_{\alpha } (\Omega _n )}{2} (y_n (s))^{-\alpha } y_n ' (s), \quad (\sigma = \varphi _n (s)).
\end{eqnarray}
{\sl Step 3: } {\sl Limit of the minimizing sequence : }
\\
Since $(x_n ' (s) )^2 + (y_n ' (s)) ^2 \equiv 1$, we obtain from (\ref{ynbound}) and  (\ref{Yn}) that  
the family $\{ Y_n ^{\alpha +1 } \} $ is equibounded and uniformly Lipschitz continuous on $(0,1)$. Hence there is a function $Y\in C[0,1] $ with $Y^{\alpha +1} \in C^{0,1} (0,1)$ such that, up to a subsequence,
\begin{equation}
\label{limYn}
Y_n ^{\alpha +1}  \longrightarrow Y^{\alpha +1 } \quad \mbox{uniformly on $[0,1]$.}
\end{equation}
Moreover, setting $Y(0)=: y^+ $, $Y(1) =: y^- $, the bounds (\ref{boundy+}) are in place and $\sigma \mapsto Y(\sigma ) $ is nonincreasing. 
\\
Let 
$$
\sigma_0 := \sup \{ \sigma\in (0, 1] : \, Y(\sigma )>0\} .
$$
Then from (\ref{estxn}) and (\ref{Xn}) we obtain that 
  the families $\{ X_n \} $ and $\{ dX_n /d\sigma \} $ are equibounded on 
every closed subset of $[0, \sigma _0)$. Hence there exists a function $X\in C[0,\sigma _0 )$ which is locally Lipschitz continuous on $[0, \sigma _0 )$, such that, up to a subsequence,
\begin{equation}
\label{limXn}
X_n \rightarrow X\quad \mbox{uniformly on closed subsets of $[0,\sigma _0 )$.}
\end{equation}
Moreover, from (\ref{estxn}) and (\ref{Xn}) we find that  
\begin{eqnarray}
\label{boundX}
\mu (\alpha , \beta ) & \geq & 2 (Y(\sigma))^{\alpha } X(\sigma ) \quad \mbox{and} 
\\
\label{boundXprime}
\mu (\alpha , \beta) +1 & \geq & 2 (Y(\sigma))^{\alpha } X^{\prime} (\sigma ), \quad   
 (\sigma \in [0,\sigma _0 )).
\end{eqnarray} 
Let $\Omega $ be the set in $\mathbb{R}^2 _+ $ with $\Omega = S(\Omega )$  such that $\partial \Omega \cap \{ (x,y):\, x>0 , y>0\} $ 
is represented by the pair of functions
$$
X(\sigma ), Y(\sigma ):\quad \sigma  \in [0, \sigma _0 )  . 
$$
In view of (\ref{boundXprime}) we have that
\begin{equation}
\label{intA}
A_{\beta } (\Omega ) = \frac{2}{\beta +1 } \int_0 ^{\sigma _0 } Y^{\beta +1} X^{\prime } \, d\sigma.
\end{equation}
\\
{\sl Step 4: } {\sl A minimizing set :}
\\
We prove that 
\begin{equation}\label{minnew}
A_\beta(\Omega)=1\,, \qquad \qquad P_\alpha(\Omega)=\mu(\alpha, \beta)\,.
\end{equation}
In order to prove the first equality, since $A_\beta(\Omega_n)=1$, we prove that
\begin{equation}\label{limnew}
\lim_{n\to +\infty}A_\beta(\Omega_n)=A_\beta(\Omega)\,.
\end{equation}
Fix some $\delta \in (0, \sigma _0 )$. 
Then (\ref{Xn}), (\ref{Yn}), (\ref{boundXprime}), (\ref{limYn}) and (\ref{limXn}) yield
\begin{eqnarray}
\label{estAbeta1}
 & &  \left| \int_0 ^{\sigma_0 -\delta } Y_n ^{\beta +1} X_n ^{\prime } \, d\sigma   - \int_0 ^{\sigma _0 -\delta } Y  ^{\beta +1} X ^{\prime} \, d\sigma \right|
\\
\nonumber
 & \leq &  \left| \int_0 ^{\sigma _0 -\delta }  Y_n ^{\beta +1 } (X_n ^{\prime} -X^{\prime})\, d\sigma \right| + \int_0 ^{\sigma _0 -\delta } |Y_n ^{\beta +1} -Y^{\beta +1}| |X^{\prime} |\, d\sigma
\\
\nonumber
 & \leq   & Y_n ^{\beta +1} |X_n -X| \Big| _0 ^{\sigma _0 - \delta } + (\beta +1) \int_0 ^{\sigma _0 -\delta } Y_n ^{\beta } |Y_n ^{\prime } | |X_n -X|\, d\sigma  
\\
\nonumber
 & & + \int_0 ^{\sigma _0 -\delta } |Y_n ^{\beta +1} -Y^{\beta +1}| |X^{\prime} |\, d\sigma 
\qquad
 \longrightarrow 0, \quad \mbox{as $n\to \infty $.} 
\end{eqnarray}
Further, (\ref{Xn}) and (\ref{boundXprime}) give
\begin{eqnarray}
\label{lim1}
\lim_{t \to 0 } \int _{\sigma _0 -t }^{\sigma _0 } (Y_n) ^{\beta +1} X_n ^{\prime } \, d\sigma & = & 0, \quad \mbox{uniformly for all $n\in \mathbb{N}$, and}
\\ 
\label{lim2}
\lim_{t \to 0 } \int _{\sigma _0 -t } ^{\sigma _0 } Y ^{\beta +1} X ^{\prime } \, d\sigma & = & 0.
\end{eqnarray}
Now (\ref{lim1}), (\ref{lim2}), (\ref{estAbeta1}) and (\ref{intA})  yield \eqref{limnew} and therefore the first of the equalities in \eqref{minnew}.

\noindent Now we prove the second inequality in  \eqref{minnew}.
With $\delta $ as above we also have
\begin{eqnarray*}
 P_{\alpha } (\Omega _n  ) & = &
2 \int_0 ^1 Y_n ^{\alpha } \sqrt{ (X_n ')^2 + (Y_n ')^2 } \, d\sigma 
\\
 & \geq & \int_0 ^{\sigma _0 -\delta }  (Y_n ^{\alpha } - Y^{\alpha }) \sqrt{ (X_n ')^2 + (Y_n ')^2 } \, d\sigma + \int_0 ^{\sigma _0 -\delta} Y ^{\alpha } \sqrt{ (X_n ')^2 + (Y_n ')^2 } \, d\sigma 
\\
 & =: & I_{n ,\delta } ^1  + I_{n, \delta } ^2 .
\end{eqnarray*} 
In view of (\ref{limYn}) and (\ref{limXn})  it follows that $\displaystyle\lim_{n\to \infty } I_{n , \delta } ^1 =0$.
\\
On the other hand, we have
$$
\liminf _{n\to \infty } I_{n, \delta } ^2  \geq  \int_0 ^{\sigma _0 -\delta}  Y ^{\alpha }  \sqrt{ (X ')^2 + (Y ')^2 } \, d\sigma .
$$ 
Define
$$
S:= \left\{ \sigma \in (0,1):\, X(\sigma )>0, \, Y(\sigma )>0 \right\} .
$$
Letting $\delta \to 0$ we obtain 
\begin{eqnarray*}
\mu (\alpha, \beta)& \geq & 
\liminf_{n\to \infty } P_{\alpha} (\Omega  _n )
\geq \int_0 ^{\sigma _0 }  Y ^{\alpha }  \sqrt{ (X ')^2 + (Y ')^2 } \, d\sigma 
\\
 & \geq & \int_S    Y ^{\alpha }  \sqrt{ (X ')^2 + (Y ')^2 } \, d\sigma 
\\
 & = & P_{\alpha } (\Omega ) \geq \mu(\alpha, \beta).
\end{eqnarray*}
Hence $\Omega $ is a minimizing set. 
\\
Note that $\Omega $ must be simply connected in view of Lemma \ref{lemma2.4}, which implies that there is a number $\sigma _1 \in [0, \sigma _0 )$ such that
\begin{equation}
\label{Sinterval}
S = (\sigma _1 , \sigma _0 ).
\end{equation}
This finishes the proof of Lemma \ref{lemma2.11}.
$\hfill \Box $
\\[0.1cm]
\hspace*{1cm}
Next we obtain differential equations for the functions $X$ and $Y$ in the proof  in Lemma \ref{lemma2.11}.

\begin{lemma}\label{lemma 2.12}Assume (\ref{cond1}) and (\ref{cond2}). Then the minimizer $\Omega $ obtained in Lemma 2.11 is bounded, and its boundary given parametrically by
\begin{eqnarray}
\label{represent}
 & & \partial \Omega \cap \{ (x,y)\in \mathbb{R}^2 :\, x> 0,\, y>0\} 
\\
\nonumber 
 & = &  \{ (x(s), y(s)):\, s\in (0,L) \} , \quad (s:\ \mbox{arclength} ),
\end{eqnarray}
where  the functions $x,y\in C^2 [0,L]$ satisfy the following  equations:
\begin{eqnarray}
\label{euler3}
 & & 
-(y^{\alpha} y')' +\alpha y^{\alpha -1} = \lambda (\beta +1) y^{\beta } x',
\\
\label{euler4}
 & &
- (y^{\alpha } x')' = - \lambda (y^{\beta +1} )' ,
\end{eqnarray} 
together with the boundary conditions 
\begin{eqnarray}
\label{bdry1}
 & & x(0) =0, \ y(0)=y^+ ,\  y'(0)=0, 
\\
\label{bdry2}
 & & y(L)= y^- , \ \mbox{ and either }
\\
\label{bdry3}
 & (i) & \ y^- >0 , \ y'(L) =0 , \ \mbox{ or}
\\
\nonumber
 & (ii) & \ y^- =0 , \ \lim_{s\nearrow L} y^{\alpha } (s) x' (s)=0,
\end{eqnarray}
for some numbers $\lambda >0$ and $0\leq y^- <y^+ $. Finally, the curve (\ref{represent}) is strictly convex. 
\end{lemma}

\noindent {\sl Proof: } We proceed in 4 steps.
\\
{\sl Step 1: } {\sl Euler equations :}
\\ 
After a rescaling the parameter $\sigma$, we see that the functions  
$X(\sigma )$ and $Y(\sigma )$ in the previous proof 
annihilate the first variation of the functional
$$
 \int_0 ^1 Y^{\alpha } \sqrt{(X') ^2 + (Y')^2 } \, d\sigma 
$$
under the constraint 
$$
 \int_0 ^1 Y^{\beta +1} X' \, d\sigma  =\mbox{ const. } (>0).
$$
Hence $X$ and $Y$ satisfy the {\sl Euler equations }
\begin{eqnarray}
\label{euler1}
 & & -\frac{d}{d \sigma } \left( \frac{Y ^{\alpha } Y' }{\sqrt{ (X') ^2 + (Y')^2 }} \right) + \alpha Y^{\alpha -1} \sqrt{ (X') ^2 + (Y')^2 } = \lambda (\beta +1) Y ^{\beta } X' ,
\\
\label{euler2} 
 & & 
-\frac{d}{d \sigma } \left( \frac{Y ^{\alpha } X' }{\sqrt{ (X') ^2 + (Y')^2 }} \right) = -\lambda \frac{d}{d\sigma } \left( Y^{\beta +1} \right),
\end{eqnarray}
where $\lambda $ is a Lagrangian multiplier, and 
\begin{equation} 
\label{posXY} 
 X(\sigma )>0, \ Y(\sigma )>0.
\end{equation}
In addition,  the following boundary conditions  are satisfied:
\begin{eqnarray}
\label{bdryXY1}
 & & X(0)=0, \ Y(0) =: y^+ >0 ,\ Y^{\prime} (0)=0, 
\\
\label{bdryXY2}
 & & \mbox{ if $Y(1)=0$ and $\lim _{\sigma \to 1} X(\sigma ) $ exists, then $\displaystyle\lim_{\sigma \to 1} Y^{\alpha } (\sigma ) X^{\prime} (\sigma ) =0$},
\\
\label{bdryXY3} 
 & & \mbox{ if $Y(1)=: y^- >0$, then $X(1)=0$ and $X'(1)=0$.}
\end{eqnarray}
{\sl Step 2: Boundedness: }
\\[0.1cm]
It will be more convenient to rewrite the above conditions in terms of the arclength parameter $s$: Set
\begin{eqnarray*}
& & s  \equiv  \psi (\sigma ):= \int_{0 } ^{\sigma}  \sqrt{ (X' (t)) ^2   + (Y'(t))^2 }\, dt ,
\\
 & & x(s)  :=  X(\psi (s)), \ y(s):= Y(\psi (s)), \ (s\in (0, L)), 
\end{eqnarray*}
where $L\in (0, +\infty ] $.
Then we have $(x'(s))^2 + (y'(s))^2 =1$ so that  (\ref{euler1}), (\ref{euler2}) yield the system of equations
 (\ref{euler3}), (\ref{euler4}). 
Integrating (\ref{euler4}) we obtain
\begin{equation}
\label{eulerextra1}
y^{\alpha } x^{\prime } = \lambda y^{\beta +1} +d,
\end{equation}
for some $d\in \mathbb{R} $.
\\
Assume first that $x(s )$ is unbounded. Then $L=+\infty $, and in view of (\ref{boundX}) we have that $\lim_{s\to \infty } y(s)=0$. 
If $\alpha =0$, then this would imply $P_0 (\Omega )=+\infty $, which is impossible. Hence we may restrict ourselves to the case $\alpha >0$. 
\\
There is a sequence $s_n \to +\infty $ such that $\displaystyle\lim_{n\to \infty } x^{\prime } (s_n )=1$. Using $s= s_n $ in  (\ref{eulerextra1}) and passing to the limit $n\to \infty $ gives $d=0$. Plugging this into  (\ref{euler3}), we find
\begin{equation}
\label{eulerextra2}
- \left( y^{\alpha } y^{\prime} \right) ^{\prime} + \alpha 
y^{\alpha -1} = \lambda ^2 (\beta +1) y^{2\beta +1 -\alpha }.
\end{equation} 
Multiplying (\ref{eulerextra2}) with $y^{\alpha } y^{\prime} $ and integrating, we  obtain
$$
- y^{2\alpha } (y^{\prime } )^2 + \alpha y^{2\alpha } = \lambda ^2 (\beta +1 ) y^{2\beta +2} + D,
$$
or equivalently,
\begin{equation}
\label{eulerextra3}
- (y^{\prime } )^2 + \alpha  = \lambda ^2 (\beta +1 ) y^{2\beta +2-2\alpha } + D y^{-2\alpha },
\end{equation}
for some $D\in \mathbb{R} $. 
Using $s=s_n $ in (\ref{eulerextra3}) and taking into account that $\lim_{n\to +\infty } y^{\prime} (s_n ) =0$, $\beta +1-\alpha >0$ and $\alpha >0$, we arrive again at a contradiction. 
Hence $x(s)$ is bounded, and we deduce the boundary conditions (\ref{bdry1})-(\ref{bdry3}). 
\\[0.1cm]
{\sl Step 3: } {\sl $\lambda$ is positive :}
\\
Multiplying (\ref{euler3}) with $y$ and integrating from $s=0$ to $s=L$ gives
$$
-\int_0^L (y^{\alpha } y')' y\, ds + \alpha \int_0 ^L y^{
\alpha } \, ds = \lambda (\beta +1) \int_0 ^L y^{\beta +1} x' \, ds .
$$
Using integration by parts this yields
$$
-y^{\alpha +1} y' \Big| _0 ^L + \int_0 ^L y^{\alpha } \left( (y') ^2 +\alpha \right) \, ds  = \lambda (\beta +1) y^{\beta +1} x \Big| _0 ^L + \lambda (\beta +1 )^2 \int_0 ^L 
y^{\beta } (-y') x \, ds .
$$
The first two boundary terms in this identity vanish due to the boundary conditions (\ref{bdry1})-(\ref{bdry3}) and the two integrals are positive since $y' \leq 0$ and $y' \not\equiv 0$. It follows that
\begin{equation}
\label{lambda>0}
\lambda >0.
\end{equation}
Now, considering
equation (\ref{eulerextra1}) at $s=L$ and taking into account the boundary conditions (\ref{bdry3}), we find
$$
0\geq  \lambda (y^-) ^{\beta +1} +d,
$$
which implies that
\begin{equation}
\label{dleq0}
d\leq 0.
\end{equation}
{\sl Step 4: } {\sl Strict convexity :}
\\
From (\ref{euler3}) and (\ref{euler4})  we obtain
\begin{eqnarray*}
\nonumber
x'' & = & \lambda (\beta +1)  y^{\beta -\alpha } y' -\frac{\alpha }{y} x'y',
\\
\nonumber
y'' & = & -\lambda (\beta +1) y^{\beta -\alpha } x' + \frac{\alpha}{y} (x')^2 .
\end{eqnarray*}
Hence, using (\ref{eulerextra1}), we find for the curvature $\kappa (s)$ of the curve $(x(s), y(s))$, ($s\in (0,L)$), 
\begin{eqnarray}
\label{curvature}
\kappa & = &
\frac{-x'y''+y'x''}{((x') ^2 + (y')^2 )^{3/2}} =   -x'y''+y'x'' 
\\
\nonumber
 & = & - \frac{\alpha}{y} x' + \lambda (\beta +1) y^{\beta - \alpha } 
  = 
\lambda y^{\beta -\alpha } ( \beta +1-\alpha ) -\alpha dy^{-1-\alpha }.
\end{eqnarray}
The last expression is positive by (\ref{lambda>0}) and (\ref{dleq0}), which means that $\Omega $ is strictly convex.   
The Lemma is proved.
$\hfill \Box $

\begin{lemma}\label{lemma2.13} Assume (\ref{cond1}) and (\ref{cond2}), and let $\partial \Omega $ be given by (\ref{represent})-(\ref{bdry3}). Then $y^- =0$.
\end{lemma}

\noindent{\sl Proof: } Supposing that $y^- >0$, we will argue by contradiction.
We proceed in 3 steps.
\\
{\sl Step 1: Another parametrization of $\partial \Omega $:}
\\
Let 
$$
x_0 := \sup \{ x: \, (x,y)\in \partial \Omega \} .
$$
Since $\Omega $ is strictly convex, there are functions $u_1 , u_2 \in C^2[0, x_0 ) \cap C[0, x_0 ] $ such that
\begin{eqnarray}
\label{domega}
 & & \partial \Omega \cap \{ (x,y): \, x>0, y>0\}
 =  \{ (x,y):\, u_1 (x)<y<u_2 (x),\, 0<x<x_0 \} ,
\\
\label{u10u20}
 & & u_1 (0) = y^- ,\ u_2 (0) =y^+ ,
\\
\label{y0x0}
 & &  u_1 (x_0 )=u_2 (x_0 ) =: y_0 ,
\\
\label{uiprime}
 & & u_1 ' (0) = u_2 '(0) =0.
\end{eqnarray}
Furthermore, the Euler equations (\ref{euler3}), (\ref{euler4}) lead to
\begin{eqnarray}
\label{euler5}
\frac{\alpha u_1 ^{\alpha -1}}{\sqrt{1+ (u_1 ')^2 }} - \frac{u _1 ^{\alpha } u_1 '' }{(1+(u_1 ')^2 )^{3/2}} & = & - \lambda (\beta +1) u_1 ^{\beta } ,
\\
\label{euler6} 
\frac{\alpha u_2 ^{\alpha -1}}{\sqrt{1+ (u_2 ')^2 }} - \frac{u _2 ^{\alpha } u_2 '' }{(1+(u_2 ')^2 )^{3/2}} & = &  \lambda (\beta +1) u_2 ^{\beta } .
\end{eqnarray}  
Using (\ref{y0x0}) and the fact that
\begin{eqnarray}
\label{inf1}
 & & \lim_{x\to x_0 } u_1 ' (x) = +\infty \quad \mbox{ and }
\\
 & &   \lim_{x\to x_0 } u_2 ' (x) = -\infty, 
\end{eqnarray}
lead to
\begin{eqnarray}
\label{euler7}
\frac{u_1 ^{\alpha }}{\sqrt{1+ (u_1 ')^2 }} & = & \lambda ( y_0 ^{\beta +1} - u_1 ^{\beta +1} ) ,
\\
\label{euler8}  
\frac{u_2 ^{\alpha }}{\sqrt{1+ (u_2 ')^2 }} & = & \lambda ( u_2 ^{\beta +1 } - y_0 ^{\beta +1}  ) .
\end{eqnarray}
Finally, the boundary conditions (\ref{u10u20})--(\ref{uiprime}) 
lead to the following formulas:
\begin{eqnarray}
\label{y0formula}
y_0 ^{\beta +1} & = & \frac{(y^-)^{\beta +1} (y^+)^{\alpha } + (y^+)^{\beta +1} (y^-)^{\alpha }}{ (y^- )^{\alpha } +  (y^+ )^{\alpha } }, 
\\
\label{lambdaformula}
\lambda & = & \frac{ (y^- )^{\alpha } +  (y^+ )^{\alpha } }{
(y^+)^{\beta +1} - (y^-)^{\beta +1}} .
\end{eqnarray}
{\sl Step 2: Curvature:} 
\\
In the following, we will refer to points $(x, u_1 (x))$ as points of the 'lower curve' and to points $(x,u_2 (x))$ as points of the 'upper curve', ($x\in (0, x_0 )$).
\\  
The signed curvature $\kappa $ (see (\ref{curvature})) can be expressed in terms of the functions $u_1 $ and $u_2 $. More precisely, we have
\begin{equation}
\label{curv2}
\kappa = 
\left\{ 
\begin{array}{ll}
 \frac{u_1 ''}{(1+ (u_1 ' )^2 )^{3/2} } & \mbox{on the lower curve}
\\
-\frac{u_2 ''}{(1+ (u_2 ' )^2 )^{3/2} } & \mbox{on the upper curve}
\end{array}
\right.
.
\end{equation} 
Accordingly, we will write
\begin{eqnarray*}
\kappa _1 (x) & := & \frac{u_1 ''}{(1+ (u_1 ' )^2 )^{3/2} }, 
\\
\kappa _2 (x) & := & -\frac{u_2 ''}{(1+ (u_2 ' )^2 )^{3/2} }, \qquad (x\in (0, x_0 )).
\end{eqnarray*}
Finally, let $s_0 \in (0,L)$ be taken such that
$y(s_0 )= y_0 $ and $x(s_0 )= x_0 $.
Then formula (\ref{eulerextra1}) taken at $s=s_0 $ leads to
\begin{equation}
\label{dform}
d= -\lambda y_0 ^{\beta +1} .
\end{equation}
Plugging this into (\ref{curvature}) we find
\begin{equation}
\label{curvaturenew}
\kappa (s) =  \lambda \left( (\beta +1-\alpha ) y^{\beta -\alpha } + \alpha y^{-1-\alpha } y_0 ^{\beta +1} \right) , \quad (s\in [0, L]).
\end{equation} 
Differentiating (\ref{curvaturenew}) we evaluate
\begin{equation}
\label{curvdiff1}
\kappa '(s) = -\lambda y^{-2-\alpha } y' \left( \alpha (\alpha +1) y_0 ^{\beta +1} - (\beta +1-\alpha ) (\beta -\alpha ) y^{\beta +1} \right)
.
\end{equation}
Since $y'(s)<0$ for $s\in (0,L)$,  
this in particular implies  
\begin{equation}
\label{curvdiff2bis}
\kappa '(s) >0 \quad \forall s\in [0,L], \qquad \mbox{
if $\beta \leq \alpha $.} 
\end{equation}
{\sl Step 4: }
\\
We claim that 
\begin{equation}
\label{uplow}
\kappa _1 (x) >\kappa _2 (x) \qquad \forall x\in (0, x_0 ).
\end{equation}
First observe that (\ref{uplow}) immediately follows from (\ref{curvdiff2bis}) if $\alpha \geq \beta $. Thus it remains to consider the case 
\begin{equation}
\label{alpha<beta}
0<\alpha <\beta < 2\alpha .
\end{equation}   
From (\ref{curvdiff1}) and the fact that
$y(s)< y_0 $ for $s\in [s_0 , L]$ we find that
\begin{eqnarray}
\nonumber 
\kappa '(s ) & \geq & -\lambda y' y^{-2-\alpha } y_0 ^{\beta +1}  \left[ \alpha (\alpha +1) - (\beta +1-\alpha )(\beta -\alpha ) \right] 
\\
\nonumber
 & = &  -\lambda y' y^{-2-\alpha } y_0 ^{\beta +1} (2\alpha -\beta ) (\beta +1)
\\ 
\label{curvdiff2}
 & > & 0 \quad \forall s\in [s_0 ,L).
\end{eqnarray}
This means that 
\begin{equation}
\label{uplow1}
\kappa _1 (x) >\kappa _2 (x) \quad \forall x\in (x_0 -\varepsilon ,x_0 ),
\end{equation}
for some (small) $\varepsilon >0$. Now assume that (\ref{uplow}) does not hold. By (\ref{uplow1}) there exists a number $x_1 \in (0,x_0 )$ such that  
\begin{eqnarray}
\label{uplow2}
 & & \kappa _1 (x) >\kappa _2 (x) \quad \forall x\in (x_1 ,x_0 ), \quad \mbox{and }
\\
\label{uplow3}
 & & \kappa_1 (x_1 )= \kappa _2 (x_1 ).
\end{eqnarray}
We claim that 
(\ref{uplow2}) implies that
\begin{equation}
\label{kappa12}
u_1 '(x_1) < -u_2 ' (x_1).
\end{equation}
To prove (\ref{kappa12}), observe first that
$$
\lim_{x\to x_0 } u_1 '(x) = \lim_{x\to x_0} (-u_2 ' (x))= +\infty .
$$ 
Then, integrating (\ref{uplow2}) over $(x_1 , x_0 )$ leads to 
\begin{eqnarray*}
 & & 
1- \frac{u_1 ' (x_1 )}{\sqrt{1+ (u_1 '(x_1 ))^2 }} 
  =  \int_{x_1 } ^{x_0 } \left( \frac{u_1 ' (x)}{\sqrt{1+ (u_1 '(x_1 ))^2 }} \right) ' \, dx 
\\
 & = & \int_{x_1 } ^{x_0 }  \frac{u_1 '' (x)}{\left[1+ (u_1 '(x_1 ))^2 \right] ^{3/2} }  \, dx = \int_{x_1 }^{x_0 }
\kappa _1 (x) \, dx 
\\
 & > & \int_{x_1 }^{x_0 }
\kappa _2 (x) \, dx 
=  \int_{x_1 } ^{x_0 }  \frac{-u_2 '' (x)}{\left[1+ (u_2 '(x_1 ))^2 \right] ^{3/2} }  \, dx 
\\
& = & \int_{x_1 } ^{x_0 } \left( \frac{-u_2 ' (x)}{\sqrt{1+ (u_2 '(x_1 ))^2 }} \right) ' \, dx
 =  1+ \frac{u_2 ' (x_1 )}{\sqrt{1+ (u_2 '(x_1 ))^2 }} ,
\end{eqnarray*}
which implies (\ref{kappa12}).
Now, from (\ref{kappa12}) we deduce  
$$
\frac{1}{\sqrt{1 + (u_1 ' (x_1 ))^2 } }
> 
\frac{1}{\sqrt{1 + (u_2 ' (x_1 ))^2 } }.
$$
Together with (\ref{euler7}) and (\ref{euler8}) we obtain from this
$$
\frac{y_0 ^{\beta +1} -(u_1 (x_1 )) ^{\beta +1}}{(u_1 (x_1 ))^{\alpha }} > 
\frac{(u_2 (x_1 )) ^{\beta +1}-y_0 ^{\beta +1} }{(u_2 (x_1 ))^{\alpha }} ,
$$
or, equivalently
\begin{equation}
\label{uplow4}
(u_1 (x_1 )) ^{\beta +1-\alpha } + 
(u_2 (x_1 )) ^{\beta +1-\alpha } < y_0 ^{\beta +1} \left( 
(u_1 (x_1 )) ^{-\alpha } + 
(u_2 (x_1 )) ^{-\alpha } \right) .
\end{equation}
Furthermore, 
multiplying (\ref{euler5}) by $(u_1 )^{-\alpha}$, respectively (\ref{euler6}) by $(u_2 )^{-\alpha }$, adding both equations  and taking into account (\ref{uplow3}) leads to
$$
\frac{\alpha}{u_1 \sqrt{1+ (u_1 ')^2 }} +
\frac{\alpha}{u_2 \sqrt{1+ (u_2 ')^2 }} = \lambda (\beta +1) \left( (u_2) ^{\beta -\alpha } - (u_1 )^{\beta -\alpha } \right) \quad \mbox{at $x=x_1 $.}
$$ 
Using once more (\ref{euler7}) and (\ref{euler8}) then gives
$$
\frac{\alpha (y_0 ^{\beta +1} -(u_1)^{\beta +1})
}{ 
(u_1 )^{\alpha +1}} +
\frac{\alpha ((u_2 )^{\beta +1} -y_0 ^{\beta +1})
}{
(u_2)^{\alpha +1}} =  (\beta +1) 
\left( 
(u_2 )^{\beta -\alpha } - (u_1 )^{\beta -\alpha } 
\right)
, 
$$
or equivalently,
$$
\alpha y_0 ^{\beta +1} 
\left( 
(u_1)^{-\alpha -1} -(u_2)^{-\alpha -1} 
\right) 
= 
(\beta +1-\alpha ) 
\left( 
(u_2 )^{\beta -\alpha } -(u_1 )^{\beta -\alpha } 
\right) \quad 
\mbox{at $x=x_1 $.}
$$
From this and (\ref{uplow4}) we then obtain
\begin{eqnarray}
\label{uplow5}
 & & (\beta +1-\alpha ) 
\left( 
(u_2 )^{\beta -\alpha } -(u_1 )^{\beta -\alpha } 
\right) 
\\
\nonumber
 & > & \alpha 
\left(
(u_1 )^{-\alpha -1} -(u_2 )^{-\alpha -1} 
\right) 
\frac{
(u_2 )^{\beta +1-\alpha } + (u_1 ) ^{\beta +1-\alpha }
}{
(u_1 )^{-\alpha } + (u_2 )^{-\alpha }
}
\quad \mbox{at $x=x_1 $.}
\end{eqnarray}
Setting
$$
z:= \frac{u_1 (x_1 )}{u_2 (x_1 )} \in (0,1),
$$
this leads to
\begin{equation}
\label{uplow6}
(\beta +1-\alpha ) 
\left(
z-z^{\beta +1-\alpha } +z^{\alpha +1} -z^{\beta +1 }
\right) 
>
\alpha
\left(
1-z^{\alpha +1} + z^{\beta +1-\alpha } -z^{\beta +2}
\right) .
\end{equation}
But this contradicts Lemma A (Appendix). This finishes the proof of (\ref{uplow}).
\\
{\sl Step 5: }
\\
Since $\lim_{x\to x_0 } u_1 ' (x) = -\lim_{x\to x_0 } u_2 ' (x) =+\infty $, (\ref{uplow}) implies that 
$$
u_1 ' (0) < -u_2 ' (0)=0,
$$
which contradicts the boundary conditions (\ref{uiprime}). 
Hence we must have that $y^- =0$.
$\hfill \Box $
\\[0.1cm]
\hspace*{1cm}Now we are in a position to give a
\\
{\sl Proof of Theorem 1.1:}
We split into two cases.
\\
{\sl Case 1:} Assume that
\begin{equation}
\label{cond22}
\beta <2\alpha .
\end{equation} 
By Lemma \ref{lemma2.13} we have $y^- =0$. 
Now, since $x'(L) =0$, equation (\ref{euler4}) at $s=L$ yields
\begin{equation}
\label{d=0}
d=0,
\end{equation}
so that
\begin{equation}
\label{odex}
x^{\prime } = \lambda y^{\beta +1 -\alpha } .
\end{equation} 
In view of Remark \ref{remark2.2} we may rescale $\Omega $ in such a way that $y^+ =1$. Then (\ref{odex}) at $s=0$ gives $\lambda =1$. Since $y^{\prime } (s)<0 $ for $s\in (0,L)$, we have $s= g(y)$ with a decreasing function $g\in C^1 (0,1)$. 
Writing 
$$
f(y) := x(g(y)),
$$
we obtain
\begin{equation}
\label{rewriteodex}
- \frac{f^{\prime } (y)}{ \sqrt{1+(f^{\prime} (y))^2 }} = y^{\beta +1-\alpha } ,
\end{equation}
and integrating this leads to \eqref{minU} and (\ref{xi}). 
\\[0.1cm]
{\sl Case (ii)} Now assume that
\begin{equation}
\label{cond23}
\beta =2\alpha .
\end{equation}
Since the case $\alpha =0$ is trivial, we may assume $\alpha >0$. Let us fix such $\alpha $.
\\
First observe that for every smooth domain $U \subset \mathbb{R}^2 _+ $, the mapping
\begin{equation}
\label{contbeta}
\beta \longmapsto {\mathcal R}_{\alpha , \beta } (U ), \quad (-1 < \beta \leq 2\alpha ),
\end{equation}
is continuous. Furthermore, from Case (i) we see that the mapping
$$
\beta \longmapsto \mu (\alpha , \beta ), \quad (-1<\beta <2\alpha ),
$$
is continuous, and the limit
$$
Z:= \lim_{\beta \to 2\alpha - } \mu (\alpha , \beta )
$$
exists. 
\\
Now let $\Omega^\star $ be the domain that is given by formulas (\ref{minU}), (\ref{xi}), with $\beta =2\alpha $. Then we also have 
$$
Z= {\mathcal R}_{\alpha , 2\alpha } (\Omega^\star ),
$$
which implies that
$ Z\geq \mu (\alpha , 2\alpha )$. 
\\
Assume that
$$
Z> \mu (\alpha , 2\alpha ).
$$
Then there is a smooth set $\Omega ' \subset \mathbb{R}^2 _+ $ such that also  
$$
Z > {\mathcal R}_{\alpha , 2\alpha }(\Omega' ) .
$$
But by (\ref{contbeta}) this implies that 
$$
{\mathcal R}_{\alpha , \beta } (\Omega ') < \mu (\alpha , \beta ),
$$
when $\beta <2\alpha $ and $|\beta -2\alpha |$ is small, which is impossible.  
Hence we have that
$$
Z= \mu (\alpha , 2\alpha )= {\mathcal R}_{\alpha , 2\alpha } (\Omega^\star ).
$$
This finishes the proof of the Theorem.
\\[0.1cm]
\hspace*{1cm}
Finally we evaluate $\mu (\alpha , \beta)$. Put $\gamma := \beta +1 -\alpha (>0)$. With the Beta function $B$ and  
the function $f$ given by (\ref{xi}) we have
\begin{eqnarray*}
P_{\alpha } (\Omega ^\star) & = & 2\int_0 ^1 y^{\alpha } \sqrt{1+ (f' (y))^2 } \, dy 
 = 
2  \int_0 ^1 
\frac{y^{\alpha } }{ 
\sqrt{1- y^{2\gamma }  } } \, dy
\\
 & = & 
\frac{1}{\gamma }
\int_0 ^1 
\frac{
z^{\frac{\alpha +1}{2\gamma } -1} \, dz
}{ 
\sqrt{1-z} 
} 
= \frac{1}{\gamma} \cdot
B
\left( 
\frac{\alpha +1}{2\gamma } , \frac{1}{2} 
\right),
\end{eqnarray*}
and
\begin{eqnarray*}
A_{\beta } (\Omega ^\star )
 & = & 
2\int_0 ^1 y^{\beta } f(y)\, dy 
= 
-\frac{2}{\beta +1 } 
\int_0 ^1 y^{\beta +1 } f ' (y)\, dy 
\\
 & = & 
\frac{2}{\beta +1} 
\int_0 ^1 
\frac{
y^{2\gamma  }
}{ 
\sqrt{1 - y ^{2\gamma }
} 
}
\, dy 
 = 
\frac{1}{\gamma (\beta +1 )} 
\int_0 ^1 
\frac{
z^{
\frac{\alpha +1}{2\gamma} }
}{ 
\sqrt{1-z}
} \, dz 
\\
 & = & 
\frac{
1
}{
\gamma (\beta +1 ) 
}
B \left( 1+ \frac{\alpha +1}{2\gamma } , \frac{1}{2} \right).
\end{eqnarray*}
Using the identity 
$$
B(a+1,b) = B(a,b) 
\cdot 
\frac{a}{a+b}, 
\qquad 
(\mbox{Re} (a)>0,\, 
\mbox{Re}(b) >0),
$$
we obtain
\begin{eqnarray*}
\mu (\alpha , \beta ) 
 & = & 
{\mathcal R}_{\alpha , \beta} (\Omega ^\star )
  =  
\frac{
P_{\alpha } (\Omega ^\star )
}{ 
\left[ 
A_{\beta } (\Omega ^\star )
\right] 
^{
\frac{\alpha +1}{\beta +2 }
}
}
\\
 & = & 
\frac{
1
}{
\gamma 
}  
\cdot 
B
\left( 
\frac{\alpha +1 }{2\gamma } , \frac{1}{2}
\right)
\cdot 
\left[
\frac{1}{\gamma (\beta )} \cdot 
B\left( 1+ \frac{\alpha +1 }{2\gamma } , \frac{1}{2}
\right)
\right] ^{
-\frac{\alpha +1}{\beta +2} 
} 
\\
 & = & 
\gamma ^{ 
\frac{\alpha +1}{\beta +2} -1 
} 
\cdot
\left[
\frac{
(\beta +1)(\beta +2)
}{
\alpha +1 
} 
\right] ^{\frac{\alpha +1}{\beta +2}} 
\cdot 
\left[ 
B 
\left( 
\frac{\alpha +1}{2\gamma } , \frac{1}{2} 
\right) 
\right] ^{\frac{\gamma}{\beta +2 }} ,
\end{eqnarray*}
which is (\ref{muab}).  In case of $\beta = 2\alpha $ this leads to 
\begin{eqnarray*}
\mu (\alpha , 2 \alpha )
 & = & 
(\alpha +1) ^{-\frac{1}{2} } 
\cdot
\left[ 2 (2\alpha +1)\right]
^{\frac{1}{2}}
\cdot  
\left[ B \left( \frac{1}{2} , \frac{1}{2} \right) 
\right] ^{\frac{1}{2} }
 =
\sqrt{\frac{2\pi (2 \alpha +1)}{ \alpha +1 } }. \qquad \qquad \Box 
\end{eqnarray*}




\begin{remark} 
\label{Remark 2.3} 
It is also well-known that the isoperimetric inequality is equivalent to the following functional inequality, (see \cite{abreu}, Lemma 3.5). 
\begin{equation}
\label{sobolev}
\iint_{\rup}|\nabla u|y^\alpha\, dxdy\ge \mu (\alpha , \beta ) \left (\iint_{\rup}| u|^\frac{\beta+2}{\alpha+1}y^\beta\, dxdy    \right )^\frac{\alpha+1}{\beta+2}\,, \quad \forall u\in C^\infty_0(\re^2)\,.
\end{equation}
\end{remark}

\section{ Applications}
In this section we firstly show that our isoperimetric inequality implies a sharp estimate of the so-called {\sl weighted Cheeger constant}. 

\noindent Then we deduce an estimate of the first eigenvalue to a degenerate elliptic Dirichlet boundary values problem. 
We begin by introducing some function spaces that will be used in the sequel.
\\
Let $\Omega $ be an open subset of $\mathbb{R}^2 _+ $ and $p\in [1, +\infty )$.
\\
By $L^p (\Omega ; y^{\beta })$ we denote the weighted 
H\"older space of measurable functions $u: \Omega \to \mathbb{R}$ such that
$$
\Vert u\Vert _{L^p (\Omega ; y^{\beta })} := 
\left( \iint_{\Omega } |u|^p y^{\beta }\, dxdy \right) ^{1/p} <+\infty .
$$
Then let $W^{1,p} (\Omega ; y^{\alpha }, y^{\beta })$ be  the weighted Sobolev space of all functions $u\in L^p (\Omega ; y^{\beta })$ possessing weak first partial derivatives which belong to $L^p (\Omega ; y^{\beta })$. A norm in $W^{1,p} (\Omega ; y^{\alpha }, y^{\beta } )$ is given by
$$
\Vert u\Vert _{W^{1,p} (\Omega ; y^{\alpha }, y^{\beta })} :=  
\Vert \, |\nabla u|\, \Vert _{L^p (\Omega ; y^{\alpha } )} + \Vert u\Vert _{L^p (\Omega; y^{\beta })}.
$$
For any function
$u\in L^1 ( \Omega ;y^{\beta })$ we write
$$
|Du| (\Omega ; y^{\alpha} ):= 
\sup 
\left\{ 
\iint_{\Omega } u \, \mbox{div}\, ({\bf v} y ^{\alpha }) \, dxdy :\,
 {\bf v }\in C^{\infty }_0 (\Omega , \mathbb{R} ^2 ),\, |{\bf v }| \leq 1 
\right\} .
$$
Then let $BV (\Omega ; y^{\alpha}, y^{\beta}) $ be the weighted BV-space of all functions $u\in L^1 (\Omega ; y^{\beta })$ such that $|Du| (\Omega ; y^{\alpha } ) < +\infty $. A norm on $BV(\Omega ; y^{\alpha}, y^{\beta}) $ is given by
$$
\Vert u\Vert _{BV(\Omega ; y^{\alpha }, y^{\beta })} := |Du| (\Omega ; y^{\alpha }) + \Vert u\Vert _{L^1 (\Omega ; y ^{\beta })}.
$$
Let us explicitly remark that for an open bounded set $\Omega \subset \mathbb{R}^2_+$ the following equality holds
$$
P_\alpha(\Omega)= |D \chi_\Omega| (\mathbb{R}^2 _+ ; y^{\alpha }) \, .
$$
Finally let $X$ be the set of all the functions $w\in C^{1}(\overline{\Omega 
})$ that vanish in a neighborhood of $\partial \Omega \cap \mathbb{R}_{+}^{2}$
. Then $V^{p}(\Omega ;y^{\alpha },y^{\beta })$ will denote the closure of $X$
in the norm of $W^{1,p}(\Omega ;y^{\alpha },y^{\beta })$.

Finally we denote by $\Omega^\bigstar$ the set $t\Omega^\star$, for $t>0$, such that $A_\beta(\Omega)=A_\beta(\Omega^\bigstar)$.

\subsection{ Weighted Cheeger sets }
We define the weighted Cheeger constant of an open bounded set $\Omega \subset \mathbb{R}^2_+$ as
\begin{equation}
\label{cc0}
h_{\alpha,\beta}(\Omega)=\inf \left \{
\frac{P_\alpha(E)}{A_\beta (E)}: \ E\subset \Omega\,,\,\, 0<A_\beta(E)<+\infty
\right\} .
\end{equation}
(see also \cite{IL}, \cite{Saracco})

We firstly prove that the existence of an admissible set which realizes the minimum in \eqref{cc0} (see also \cite{Saracco}).

\begin{lemma}\label{lemma3.1} Assume $0\le \alpha<\beta+1$ and $\beta \le 2\alpha$. For any open bounded set $\Omega \subset\subset \mathbb{R}^2_+$, there exists at least one set $M\subseteq \Omega$, the so-called {\sl weighted Cheeger set}, such that 
\begin{equation}\label{cs}
h_{\alpha,\beta}(\Omega)=\frac{P_\alpha(M)}{A_\beta (M)}\,.
\end{equation}
\end{lemma}

\noindent{\sl Proof: } Since $\Omega$ is open, $h_{\alpha,\beta}(\Omega)$ is finite: Indeed, it is easy to verify that for any ball  $B$ with $B\subset \subset \Omega$, the ratio $\frac{P_\alpha(B)}{A_\beta (B)}$ is finite.

\noindent Let $\{E_k\}$ be a minimizing sequence for \eqref{cc0}. Since $\Omega$ is bounded, we have
$$
A_\beta(E_k)=\iint\limits_{E_k  }y^{\beta }\, dxdy\le A_\beta(\Omega)=\iint\limits_{\Omega  } y^{\beta }\, dxdy<+\infty\,.
$$ 
Now fix $\varepsilon>0$. There exists an index $\overline k$ such that
$$
\left | h_{\alpha,\beta}(\Omega)- \frac{P_\alpha(E_k )}{A_\beta (E_k )}  \right |<\varepsilon\, , \qquad \forall k>\overline k.
$$ 
Since $\Omega$ is bounded, for all $k>\overline k$,  we get
$$
P_\alpha(E_k )< (\varepsilon+h_{\alpha,\beta}(\Omega))A_\beta (E_k )\le (\varepsilon+h_{\alpha,\beta}(\Omega))A_\beta (\Omega )\equiv C
$$
This implies
$$
|D\chi_{E_k}(\mathbb{R}^2 _+, y^\alpha)|=P_\alpha(E_k)\le C\, \, \qquad \forall k>\overline k.
$$
Hence $\{E_k\}$ is an equibounded family in weighted $BV(\Omega; y^\alpha, y^\beta)-$norm. Thus by Lemma B (Appendix), up to subsequences, $\{\chi(E_k)\}$ converges in the weighted $L^1(\Omega; y^\beta)-$norm and pointwise a.e. to a function $u$. Moreover there exists a subset $M\subseteq \Omega$ such that $u=\chi_{M}$.

\noindent Since $\{E_k\}$ is a minimizing sequence, by lower semicontinuity of perimeter $P_\alpha$ and Lebesgue dominated convergence theorem, we get
\begin{equation}\label{limit}
h_{\alpha,\beta}(\Omega)=\lim_{k\to +\infty}\frac{P_\alpha(E_k)}{A_\beta (E_k)}\ge \frac{P_\alpha(M)}{A_\beta (M)} \,.
\end{equation}
It remains to prove that $M$ is an admissible set, that is we need to prove that 
\begin{equation}\label{pos}
A_\beta (M)>0.
\end{equation}
Assume by contradiction that $A_\beta (M)=0$. This implies that $\displaystyle\lim_{k\to +\infty}A_\beta (E_k)=0$. 

\noindent Now for a fixed $\eta>0$ consider the set 
$$
E_{k, \eta}=\{(x,y)\in E_k\,:\,\, y>\eta  \}\,.
$$
Now the following inequality holds true
$$
A_\beta (E_{k, \eta})< \delta^\beta A_0 (E_{k, \eta})\,, \quad 
P_\alpha(E_{k, \eta})< \eta^\alpha P_0(E_{k, \eta})\,,
$$
where $\delta=\eta$, if $\beta\le 0$ and $\delta=R$ for  a suitable $R>0$ such that a ball $B_R$ of radius $R$ contains $\Omega$,  if $\beta> 0$.
Denote by $B_{r_k,\eta}$  a ball of radius $r_k$ having the same Lebesgue measure $A_0 (E_{k, \eta})$ of $E_{k, \eta}$. By the classical isoperimetric inequality, we get
$$
\frac{P_\alpha(E_{k, \eta})}{A_\beta (E_{k, \eta})}
\ge \frac{ \eta^\alpha P_\alpha(E_{k, \eta})}{\delta^\beta A_0 (B_{r_k,\eta})}
 = \frac{\eta^\alpha P_\alpha(B_{r_k,\eta})}{\delta^\beta A_0 (B_{r_k,\eta})}=\frac{2 \eta^\alpha}{\delta^\beta r_k}\rightarrow +\infty\,,\quad \hbox{as $k$ goes to $+\infty$}\,.
$$
This yields a contradiction. Therefore \eqref{pos} holds true and the conclusion follows.
$\hfill \Box $
\\[0.1cm]
Once we have proved the existence of a weighted Cheeger set, we can obtain the following result.
\\[0.1cm]
\begin{theorem}\label{theorem3.2}Assume $0\le \alpha<\beta+1$ and $\beta \le 2\alpha$ and $\alpha< \beta+1$, and let $\Omega$ be a  bounded  open subset $\mathbb{R}^2_+$. Then 
 the following estimate holds true
\begin{equation}\label{ccth}
h_{\alpha, \beta}(\Omega)\ge h_{\alpha, \beta}(\Omega ^\bigstar)=
\frac{P_\alpha(\Omega ^\bigstar)}{A_\beta (\Omega ^\bigstar)}\,.
\end{equation}
\end{theorem}

\noindent {\sl Proof: }  Let $E$ be a nonempty subset of $\Omega$ with $P_\alpha(E)<+\infty$. 
By our isoperimetric inequality Theorem \ref{maintheorem} and since $E^\star\subset \Omega^\bigstar$ with $A_\beta (E)=A_\beta (E^\bigstar)$, we have that  
\begin{equation}\label{cc}
\frac{P_\alpha(E)}{A_\beta (E)}\ge \frac{P_\alpha(E^{\bigstar })}{A_\beta (E^{\bigstar })}\ge h_{\alpha,\beta}(\Omega^\bigstar)\,.
\end{equation}
It remains to prove the equality in \eqref{ccth}. Let $F$ be a nonempty subset of $\Omega^\star$. Then we have $F^\bigstar \subset \Omega^\bigstar$ and 
\begin{equation}\label{cc1}
\frac{P_\alpha(F)}{A_\beta (F)}\ge \frac{P_\alpha(F^{\bigstar })}{A_\beta (F^{\bigstar })}=t^{\beta+1-\alpha}\frac{P_\alpha(tF^{\star })}{A_\beta (tF^{\star })}
\end{equation}
for all $t>0$.
Since $F^\bigstar \subset \Omega^\bigstar$, there exists $t\ge 1$ such that $tF^\star=\Omega^\bigstar$.
Therefore 
$$
\frac{P_\alpha(F)}{A_\beta (F)}\ge t^{\beta+1-\alpha}\frac{P_\alpha(\Omega^{\bigstar })}{A_\beta (\Omega^{\bigstar })}
$$
which proves the equality in \eqref{ccth}.
$\hfill \Box $
\\[0.1cm]
%
%
{\bf Remark 3.1.}  
 Theorem 3.1 could be stated as an estimate of   the first eigenvalue of the weighted 1-laplacian.

\subsection{ A nonlinear eigenvalue problem}

\noindent Let $\Omega \subset \mathbb{R}_{+}^{2}$ be a bounded domain and
let $p\in (1,+\infty ).$ We consider the following weighted eigenvalue
problem
\begin{equation}
\left\{ 
\begin{array}{ccc}
-\text{div}\left( y^{p\gamma _{1}}\left\vert \nabla u\right\vert
^{p-2}\nabla u\right) =\lambda y^{\frac{p\gamma _{2}}{p-1}}\left\vert
u\right\vert ^{p-2}u & \text{in} & \Omega  \\ 
&  &  \\ 
\hspace{1.7 cm} u=0 & \text{on} & \partial \Omega \cap \mathbb{R}_{+}^{2},
\end{array}
\right.   \label{P}
\end{equation}
where 
\begin{equation*}
\gamma _{1}=\alpha -\frac{p-1}{p}\beta \text{ \ and \ }\gamma _{2}=\frac{p-1
}{p}\beta .
\end{equation*}
By a solution to problem \eqref{P} we mean a function \thinspace $u\in
V^{p}\left( \Omega ;y^{p\gamma _{1}}, y^{\frac{p\gamma _{2}}{p-1}}\right) $
such that
\begin{equation*}
\iint\limits_{\Omega  }
\left\vert \nabla u\right\vert ^{p-2}\nabla u\nabla \psi
y^{p\gamma _{1}}dxdy=\lambda \iint\limits_{\Omega  }\left\vert u\right\vert
^{p-2}u\psi y^{\frac{p\gamma _{2}}{p-1}}dxdy
\end{equation*}
for all function $\psi \in C^{1}(\overline{\Omega })$ such that $\psi =0$ on 
$\partial \Omega \cap \mathbb{R}_{+}^{2}.$

Let us denote by $T$ the range of values of $\alpha$
and $ \beta $ for which the isoperimetric inequality holds true. We have that 
\begin{equation}
(\alpha ,\beta )\in T
=
\left\{ \alpha \geq 0\right\}
\cap
 \left\{ \beta >\alpha -1\right\} \cap \left\{
\beta \leq  2\alpha \right\}  
\label{T_ab}
\end{equation}
if and only if
\begin{equation}
(\gamma _{1},\gamma _{2})\in 
U
:=
\left\{ \gamma_{1}+\gamma _{2} \geq 0 \right\} 
\cap 
\left\{ \frac{p  \gamma _{2}}{p-1} > \gamma _{1}+\gamma _{2}-1 \right\}   
\cap 
\left\{  \frac{p  \gamma _{2}}{p-1} \leq 2(\gamma _{1}+\gamma _{2}) \right\} .  \label{T_12}
\end{equation}

Furthermore the smallest eigenvalue of problem \eqref{P}, $\lambda
_{1,p}^{\gamma _{1},\gamma _{2}}(\Omega ),$ has the following variational
characterization 
\begin{equation*}
\lambda _{1,p}^{\gamma _{1},\gamma _{2}}(\Omega )=\min \left\{ \frac{
\displaystyle\displaystyle\iint\limits_{\Omega }\left\vert \nabla
u\right\vert ^{p}y^{p\gamma _{1}}dxdy}{\displaystyle\displaystyle
\iint\limits_{\Omega }u^{p}y^{\frac{p\gamma _{2}}{p-1}}dxdy}\text{ with }u\in
V^{p}\left( \Omega ;y^{p\gamma _{1}}, y^{\frac{p\gamma _{2}}{p-1}
}\right) \backslash \left\{ 0\right\} \right\} .
\end{equation*}
Indeed,  see e.g. Theorem 8.9 in \cite{GO}, 
for any $(\gamma _{1},\gamma _{2})\in U,$ 
  the following compact 
embedding holds true
\begin{equation*}
V^{p}\left( \Omega ;y^{p\gamma _{1}},y^{\frac{p\gamma _{2}}{p-1}
} \right) \hookrightarrow \hookrightarrow L^{p}\left( \Omega ;y^{\frac{
p\gamma _{2}}{p-1}} \right) .
\end{equation*}

\bigskip
By adapting the arguments used in \cite{Ch}, \cite{KF}, we obtain the following result

\begin{theorem}  \label{PDEs}
Let $(\gamma _{1},\gamma _{2})\in
U,$ then the following estimate holds true
\begin{equation*}
\lambda _{1,p}^{\gamma _{1},\gamma _{2}}(\Omega )\geq 
\frac{1}{p^{p}}\left[ h_{\alpha,\beta}(\Omega^\bigstar)\right] ^{p}=
\frac{1}{p^{p}}\left[ 
\frac{P_{\alpha }( \Omega ^{\bigstar }) }{A_{\beta }( \Omega
^{\bigstar }) }\right] ^{p}\,.
\end{equation*}
\end{theorem}

\noindent
\\
{\bf Remark 3.2.} 
 If $\alpha =\beta >0$ and $p=2$, a
Faber-Krahn type inequality for $\lambda _{1,p}^{\gamma _{1},\gamma
_{2}}(\Omega )$ holds true (see \cite{MadSalsa}).
Indeed in this case, we have $\alpha =2\gamma _{1}=2\gamma _{2}=\beta.$

\bigskip
\noindent {\sl Proof of Theorem \ref{PDEs}  }
\noindent We claim that
\begin{equation}
\label{claim}
\lambda _{1,p}^{\gamma _{1},\gamma _{2}}(\Omega )\geq \frac{\left[ h_{\alpha
,\beta }\left( \Omega \right) \right] ^{p}}{p^{p}},
\end{equation}
where $ h_{\alpha
,\beta }\left( \Omega \right)$ is defined in \eqref{cc0}.
Let $u$ be an eigenfunction corresponding to $\lambda _{1,p}^{\gamma
_{1},\gamma _{2}}(\Omega )$. H\"{o}lder's inequality gives
\begin{eqnarray*}
\displaystyle\iint\limits_{\Omega }\left\vert \nabla u\right\vert
u^{p-1}y^{\alpha }dxdy 
&=&
\displaystyle\iint\limits_{\Omega }\left\vert
\nabla u\right\vert u^{p-1}y^{\gamma _{1}}y^{\gamma _{2}}dxdy \\
&\leq &
\left( \displaystyle\iint\limits_{\Omega }\left\vert \nabla
u\right\vert ^{p}y^{p\gamma _{1}}dxdy\right) ^{\frac{1}{p}}\left( 
\displaystyle\iint\limits_{\Omega }u^{p}y^{\frac{\gamma _{2}p}{p-1}
}dxdy\right) ^{\frac{p-1}{p}},
\end{eqnarray*}
and therefore
\begin{equation*}
\displaystyle\iint\limits_{\Omega }\left\vert \nabla u\right\vert
^{p}y^{p\gamma _{1}}dx\geq \frac{\left( \displaystyle\iint\limits_{\Omega
}\left\vert \nabla u\right\vert u^{p-1}y^{\alpha }dxdy\right) ^{p}}{\left( 
\displaystyle\iint\limits_{\Omega }u^{p}y^{\beta }dxdy\right) ^{p-1}}
=
\frac{\left( 
\displaystyle\iint\limits_{\Omega }\left\vert \nabla \left( u^{p}\right)
\right\vert y^{\alpha }dxdy\right) ^{p}}{p^{p}\left( \displaystyle
\iint\limits_{\Omega }u^{p}y^{\beta }dxdy\right) ^{p-1}}.
\end{equation*}
If we set $f:=u^{p}$, then the previous inequality becomes
\begin{equation*}
\frac{\displaystyle\iint\limits_{\Omega }\left\vert \nabla u\right\vert
^{p}y^{p\gamma _{1}}dxdy}{\displaystyle\iint\limits_{\Omega }u^{p}y^{\beta
}dxdy}\geq \frac{\left( \displaystyle\iint\limits_{\Omega }\left\vert \nabla
f\right\vert y^{\alpha }dxdy\right) ^{p}}{p^{p}\left( \displaystyle
\iint\limits_{\Omega }fy^{\beta }dxdy\right) ^{p}}.
\end{equation*}
On the other hand Coarea formula yields 

\begin{equation*}
\frac{\left( \displaystyle\iint\limits_{\Omega }\left\vert \nabla
f\right\vert y^{\alpha }dxdy\right) ^{p}}{p^{p}\left( \displaystyle
\iint\limits_{\Omega }fy^{\beta }dxdy\right) ^{p}}=\frac{1}{p^{p}\left( 
\displaystyle\iint\limits_{\Omega }fy^{\beta }dxdy\right) ^{p}}\left(
\int_{0}^{\max f}P_{\alpha }\left( \left\{ f(x)>t\right\} \right) dt\right)
^{p}
\end{equation*}
\begin{equation*}
\geq \frac{1}{p^{p}\left( \displaystyle\iint\limits_{\Omega } 
f y^{\beta }dxdy\right) ^{p}}\left(
\int_{0}^{\max f}\frac{P_{\alpha }\left( \left\{ f(x)>t\right\} \right) dt}{
A_{\beta }\left( \left\{ f(x)>t\right\} \right) }A_{\beta }\left( \left\{
f(x)>t\right\} \right) dt\right) ^{p}
\end{equation*}
\begin{equation*}
\geq \frac{1}{p^{p}\left(\displaystyle\iint\limits_{\Omega } fy^{\beta }dxdy\right) ^{p}}\left(
\int_{0}^{\max f}h_{\alpha ,\beta }\left( \left\{ f(x)>t\right\} \right)
A_{\beta }\left( \left\{ f(x)>t\right\} \right) dt\right) ^{p}
\end{equation*}
\begin{equation*}
\geq \frac{\left[ h_{\alpha ,\beta }\left( \Omega \right) \right] ^{p}}{
p^{p}\left( \displaystyle\iint\limits_{\Omega } f y^{\beta }dxdy\right) ^{p}}
\left( \int_{0}^{\max f}A_{\beta }\left( \left\{ f(x)>t\right\} \right)
dt\right) ^{p}
=
\frac{\left[ h_{\alpha ,\beta }\left( \Omega \right) \right]
^{p}}{p^{p}}.
\end{equation*}

The claim is hence  proved. It immediately implies
 Theorem \ref{PDEs}, thanks to  \eqref{ccth}.

\bigskip

\section{Appendix} 
We prove two technical results.
\\[0.1cm]
{\bf Lemma A:} {\sl Let $\alpha , \beta \in \mathbb{R} $ with $0<\alpha \leq\beta \leq 2\alpha $.
Then}
\begin{eqnarray}
\nonumber & &
\alpha (1-z^{\alpha +1} + z^{\beta +1-\alpha }-z^{\beta +2} ) + (\beta +1 -\alpha ) (-z +z^{\beta +1-\alpha } -z^{\alpha +1} +z^{\beta +1} )
\\
\label{app1}
 &  &  \quad > 0 \qquad \forall z\in (0,1).
\end{eqnarray}
{\sl Proof: } We fix $\alpha >0$ and define
\begin{eqnarray*}
g(z,\beta ) & := & \alpha (1-z^{\alpha +1} + z^{\beta +1-\alpha }-z^{\beta +2} ) + (\beta +1 -\alpha ) (-z +z^{\beta +1-\alpha } -z^{\alpha +1} +z^{\beta +1} ), 
\\
 & & \mbox{and} 
\\
 h(z) & := & g(z, 2\alpha ) = \alpha (1-z^{2\alpha +2}) + (\alpha +1 ) (-z + z^{2\alpha +1}), 
\\
 & &   (z\in [0,1],\ \beta \in [\alpha , 2 \alpha ] ).
\end{eqnarray*}
Then 
\begin{eqnarray}
\label{app2}
 & & g(z, \alpha )  =  \alpha (1-z^{\alpha +1} )(1+z) >0,
\\
\label{app3}
 & & h(0)=\alpha , \ h(0)=0,
\\
\nonumber 
 & & h'(z) = -2\alpha (\alpha +1) z^{2\alpha +1} + (\alpha +1) [ -1 +( 2\alpha +1) z^{2\alpha } ] \quad \mbox{and }
\\
\label{app4}
 & & h'' (z) = 2\alpha (\alpha +1) (2\alpha +1) z^{2\alpha +1} (1-z) >0.
\end{eqnarray}
Hence $h'(1)=0$ which together with (\ref{app3}) and (\ref{app4}) implies that 
\begin{equation}
\label{app5}
h(z) = g(z, \alpha ) >0 \quad \forall z\in (0,1).
\end{equation}
Furthermore we have
\begin{eqnarray*}
\frac{\partial g}{\partial \beta } (z, \beta ) & = & z^{\beta +1-\alpha } \left[ \alpha (1-z^{1+\alpha }) + (\beta +1-\alpha ) (1+z^{\alpha }) \right]  
 \ln z 
\\
 & & -(1-z^{\beta }) (z+z^{\alpha +1}) 
\\
 &  < & 0 \qquad \mbox{if $z\in (0,1) $ and $\beta \in [\alpha , 2 \beta ] $.} 
\end{eqnarray*}
Together with (\ref{app2}) and (\ref{app4}) this implies
$$
g(z,\alpha )>0 \quad \mbox{if $z\in (0,1)$ and $\beta \in [\alpha , 2 \alpha ],$}
$$
which is (\ref{app1}).  
$\hfill \Box $
\\[0.1cm]
{\bf Lemma B: }{\sl Let $\{u_n\}\subset BV(\Omega; y^\alpha, y^\beta)$ be a bounded sequence. Then there exists a subsequence that converges in $L^1(\Omega; y^\beta)$ and a.e. in $\Omega$ to some function $u$.}
\\[0.1cm]
{\sl Proof: } Put $\gamma =\frac{\beta+2}{\alpha+1}(>1)$ and let $\Omega_\varepsilon=\Omega \cap \{ (x,y)\,:\, y>\varepsilon \}$ for any $\varepsilon >0$. 
\\
Let $k\in \N$.
By a classical compactness result in the unweighted case, there exists a function $u^k\in L^1(\Omega_{2^{-k}}; y^\beta)$ and an increasing sequence of integers $\{ a(k,m)\}_{m\geq 1 } $ such that 
\begin{equation}
\label{subseq}
u_{a(k,m)}\rightarrow u^k \qquad \hbox{in }L^1(\Omega_{2^{-k}}; y^\beta) \hbox{ and a.e. in }\Omega\,.
\end{equation}
By choosing $\{ a(k+1,m)\}$ to be a subsequence of $\{ a(k,m\} $, ($k\in \N$), we can achieve that $u^k=u^{k+1}$ in $\Omega_{2^{-k}}$, $k\in \N$.
\\ 
Now put 
$$
u(x) = \left\{
\begin{array}{ll}
u^1(x) 
 & 
\mbox{if  } x\in \Omega_{2^{-1}}
\\
u^k (x) 
 & 
\mbox{if  } x\in 
\Omega_{2^{-k}}\setminus \overline\Omega_{2^{-k+1}}\,
, k=2,3,....
\end{array}
\right.
$$
In view of our isoperimetric inequality, the sequence $u_n $ is equibounded in $L^\gamma (\Omega; y^\beta)$.
We have the following estimate:
\begin{eqnarray*}
\| u \| _{L^1(\Omega_{2^{-k}}; y^\beta)}
 & = & \sum_{j=0}^k
\int_{
\Omega_{2^{-j}}
 \setminus \overline\Omega_{2^{-j+1}}
 }
|u| y^\beta\, dxdy\\
 & \le & 
 \sum_{j=0}^k
\left( \int_{
\Omega_{2^{-j}}\setminus \overline\Omega_{2^{-j+1}}
}|u|^\gamma y^\beta
\, dxdy\right )^\frac 1\gamma
\left ( \int_{\Omega_{2^{-j}}\setminus \overline\Omega_{2^{-j+1}}}y^\beta\, dxdy\right )^{1-\frac 1\gamma}\\
 & \le &
\sum_{j=0}^k C_1 \left [  y^{\beta+1}\Bigg |^{2^{-j+1}}_{2^{-j}}  \right ]^{
1-\frac 1\gamma}
\le C_2 \sum_{j=0}^\infty 2^{(-j+1)(\beta+1)(1-\frac 1\gamma)}< +\infty \, ,
\end{eqnarray*}
with constants that do not depend on $k$. This implies that $u\in L^1(\Omega; y^\beta)$.
\\
Let $\varepsilon >0$. Choose $k$ large enough such that 
$$
\int_{ \Omega\setminus \overline\Omega_{2^{-k}}
} y^\beta dxdy
 < 
\varepsilon
$$ 
and then $m$ large enough such that
$$
\int_{
\Omega_{2^{-k}}
}|u-u_{a(k,m)} | y^\beta dxdy
<
\varepsilon
\, .
$$ 
Then we obtain
\begin{eqnarray*}
\int_{\Omega}|u-u_{a(k,m)} | y^\beta dxdy
 & < &
\varepsilon+
\int_{
\Omega \setminus\overline\Omega_{2^{-k}}
} 
|u-u_{a(k,m)}| y^\beta dxdy
\\
 & \le & 
\varepsilon + \left(  \int_\Omega |u-u_{a(k,m)}| y^\beta dxdy    \right)
^\frac 1\gamma
\left(  \int_{
\Omega \setminus\overline\Omega_{2^{-k}}
}  
y^\beta dxdy    \right)^{1-\frac 1\gamma}\\
 & \le & 
\varepsilon +C \varepsilon ^{(1-\frac 1\gamma)(\beta+1)}\,
,
\end{eqnarray*}
where $C$ does not depend on $k$. From this 
the assertion follows. 
$\hfill \Box $

\vspace{0.2 cm}

{\bf Acknowledgement: } We are grateful to M. van den Berg  for giving us useful suggestions. The first, third, fourth and fifth authors are members of the Gruppo Nazionale per l'Analisi Matematica, la Probabilit\`a
e le loro Applicazioni (GNAMPA) of the Istituto Nazionale di Alta Matematica (INdAM) and they thank this institution for the support. The second author was
partially supported by Leverhulme Trust ref. VP1-2017-004. We are also grateful to the Departments of Mathematics of Swansea University and of the University of Naples Federico II, and to South China University of Technology (ISAM) at Guangzhou for visiting appointments and their colleagues for their
kind hospitality.

\end{document}